\documentclass[12pt]{amsart}
\usepackage{amssymb}
\usepackage[all]{xy}
\usepackage{graphicx}

\providecommand{\noglossaryignore}[1]{}
\newcommand{\globalglossaryentry}[3]{\makebox[1.5in][l]{\tt $\backslash${#1}} 
\makebox[1.1in][l]{{$#2$}} \makebox[2.5in][l]{{#3}}\newline} 
\newcommand{\newcommandabbreviation}[3]{\newcommand{#1}{#2}%
\noglossaryignore{\globalglossaryentry{#3}{#2}{}}}
\newcommand{\renewcommandabbreviation}[3]{\renewcommand{#1}{#2}%
\noglossaryignore{\globalglossaryentry{#3}{#2}{}}}
\newcommand{\newcommandmacro}[4]{\newcommand{#1}{#2}%
\noglossaryignore{\globalglossaryentry{#3}{#2}{#4}}}
\newcommand{\gge}[3]{\noglossaryignore{\globalglossaryentry{#1}{#2}{#3}}}
\newcommand{\myaddress}%
{\parbox{3in}{\footnotesize \begin{center} 
Mathematics Department, City University, \\  
Northampton Square, London EC1V 0HB, UK.\end{center}}}
\newcounter{minidef}[section]
\renewcommand{\theminidef}{\thesection.\arabic{minidef}}
\newcommand{\mdef}{\refstepcounter{minidef} 
\medskip \noindent ({\bf \theminidef}) }
\newcounter{minicapt}

\newtheorem{de}{Definition}     \newtheorem{pr}{Proposition} 
   
\newtheorem{lem}{Lemma}

\noglossaryignore{GREEK ETC.\newline}
\newcommandabbreviation{\e}{\epsilon}{e}        
\newcommandabbreviation{\lam}{\lambda}{lam}  
\newcommandabbreviation{\la}{\langle}{la}        
\newcommandabbreviation{\ran}{\rangle}{ran}
\newcommandabbreviation{\ha}{\#}{ha}             
\newcommandabbreviation{\rmap}{\rightarrow}{rmap}
\newcommandabbreviation{\aaa}{\alpha}{aaa}        
\newcommandabbreviation{\ab}{\alpha,\beta}{ab}
\newcommandabbreviation{\aab}{a(\ab )}{aab}       
\noglossaryignore{\newline RINGS\newline}
\newcommandabbreviation{\HH}{H \!\!\! I}{HH}               % Hecke algebra
\newcommandabbreviation{\C}{\mathbb C}{C}
\newcommandabbreviation{\N}{\mathbb N}{N}   %new AMS versions (was \Bbb)!
\newcommandabbreviation{\Z}{\mathbb Z}{Z}      % AMS versions!!!!!!
\renewcommandabbreviation{\Re}{\mathbb R}{Re}
\newcommandabbreviation{\R}{{\mathbb R}}{R}
\newcommandabbreviation{\Q}{\mathbb Q }{Q}
\renewcommandabbreviation{\H}{\mathbb H }{H}
\noglossaryignore{\newline SYMMETRIC GROUP\newline}
\def\Sym(#1){\Sigma(#1)}                   % generic symbol for symmetric group
\gge{Sym(-)}{\Sym(-)}{symmetric group on - objects}
\def\Sy(#1){\Sigma_{#1}}                   % symmetric group irrep.
\gge{Sy(-)}{\Sy(-)}{symmetric group irreducible -}
\def\sym(#1){\mbox{\LARGE s}(#1)}        % another generic symbol for
\gge{sym(-)}{\sym(-)}{symmetric group on - objects (variant)}
\def\sy(#1){\mbox{\LARGE s}({#1})}        % another symm gp irrep. 
\gge{sy(-)}{\sy(-)}{symmetric group irreducible - (variant)}
\newcommandmacro{\cs}{\C \, \sy(n)}{cs}{symmetric group algebra over $\C$}
\noglossaryignore{\newline PARTITIONS/SETS\newline}
\newcommand{\Nset}[1]{\underline{#1}}
\gge{Nset\{-\}}{\Nset{-}}{set of natural numbers to -} %%% (default n)}
\def\nset(#1){ \{ #1 \}_{ \underline{n} }} % the set {#1}_{n}
\gge{nset(-)}{\nset(-)}{a set $-\times\Nset$}
\def\ul(#1){_{\underline{#1}}}             % _underline #1
\gge{ul(-)}{{}\ul(-)}{subscript underline -}
\def\Ee(#1){{\bf E}_{#1}}                  % set of equiv. relations
\gge{Ee(-)}{\Ee(-)}{set of equivalence relations on set -}
\def\Eee(#1){{\bf E}_{\{ #1 \}_{\underline{n}}}}   %ditto for nset
\gge{Eee(-)}{\Eee(-)}{ditto for nset}
\def\Een(#1,#2){{\bf E}_{\{ #1 \}_{\underline{#2}}}}   %ditto for n+1set
\def\Ssn(#1,#2){{\bf S}_{\{ #1 \}_{\underline{#2}}}}   %partitions for n+1set
\def\Ss(#1){{\bf S}_{#1}}                  % set of partitions
\def\Sss(#1){{\bf S}_{\{ #1 \}_{\underline{n}}}}   %ditto for nset
\def\bbc(#1){((\beta_1)(\beta_2)...(\beta_{#1}))}      % beta singletons 
\newcommandmacro{\Ln}{{\Gamma}^{n}}{Ln}{large index set}
\newcommandmacro{\LnQ}{{\Gamma}^{n}_Q}{LnQ}{index set}
\newcommandmacro{\Zz}{\zeta}{Zz}{`shape' function}
\noglossaryignore{\newline PARTITION ALGEBRA\newline}
\def\ka(#1){\kappa_{#1}}                   % maps AP_{n=#1}A to P_{n-1}
\def\Sm(#1){\Sigma_{#1}}                   % image of S_lamda in P_n
\newcommandmacro{\com}{\bullet}{com}{bullet composition}
\newcommandmacro{\enm}{\; e^n(\! m\! ) \;}{enm}{product of idempotents}
\def\Ai(#1){ A^{ #1 \cdot } }              % A_i
\def\Aij(#1,#2){ A^{ #1  #2 } }            % A_ij
\newcommandmacro{\One}{\mbox{\bf $1 \!\!\! 1$}}{One}{algebra unit 1}
\newcommandmacro{\Bp}{B_p}{Bp}{partition basis}
\def\Bb(#1){B_p[#1]}                       % partition basis
\def\Pp(#1){P_n[#1]}                       % left module
\def\Ps(#1){P_n[#1] \! /}                  % left module
\newcommandmacro{\Ph}{\hat{P}}{Ph}{P hat  algebra}
\def\Is(#1){\sim^{#1}}                     % is equivalent under a to
\noglossaryignore{\newline MODULES\newline}
\def\Wm(#1){{\cal S}_{#1}}                 % Weyl module S 
\gge{Wm(-)}{\Wm(-)}{Weyl module with index -}
\def\wm(#1,#2){{}_{#1}{\cal S}_{#2}}       % Weyl module nS
\gge{wm(-1,-)}{\wm(-1,-)}{Weyl module with index -}
\def\Ind(#1,#2,#3){\mbox{Ind}_{#1}^{#2}#3} % induction
\gge{Ind(-1,-2,-)}{\Ind(-1,-2,-)}{induction}
\def\Res(#1,#2,#3){\mbox{Res}_{#1}^{#2}#3} % restriction
\gge{Res(-1,-2,-)}{\Res(-1,-2,-)}{restriction}
\newcommandabbreviation{\weyl}{standard}{weyl}
\newcommandabbreviation{\head}{\mbox{head }}{head}
\newcommandabbreviation{\Weyl}{Weyl}{Weyl}
\def\SS(#1){{\cal S}_{#1}}                 % Specht/Weyl module
\gge{SS(-)}{\SS(-)}{Specht/Weyl module index -}
\def\LL(#1){{\cal L}_{#1}}                 % Simple module
\gge{LL(-)}{\LL(-)}{simple module index -}
\noglossaryignore{\newline FUNCTORS/MAPS\newline}
\newcommandmacro{\Gg}{{\cal G}}{Gg}{G Functor}
\newcommandmacro{\Fg}{{\cal F}}{Fg}{F Functor}
\newcommandmacro{\ra}{\rightarrow}{ra}{}
\def\ses(#1,#2,#3){0\ra #1 \ra #2 \ra #3 \ra 0}   %short exact sequence
\gge{ses(1,2,-)}{\ses(1,2,-)}{\hspace{.5in} short exact sequence}
\def\starr(#1){ \stackrel{ #1 }{\longrightarrow} }
\gge{starr(-)}{\starr(-)}{}
\newcommandmacro{\doublerightarrow}{\; -\!\!\! -\!\!\!\!\!\! \gg \;}
{doublerightarrow}{}%{ $---->>$ }
\noglossaryignore{\newline PARTITION ALGEBRA MAPS\newline}
\newcommandmacro{\smap}{s}{smap}{`inclusion' map}
\newcommandmacro{\tmap}{t}{tmap}{$ P_n -> S_n$}
\newcommandmacro{\pmap}{\psi}{pmap}{$ S_n -> P_n $}
\noglossaryignore{\newline MISC.\newline}
\def\Amap(#1){{\cal A}_{#1}}               % variant inclusion A_Gamma
\gge{Amap(-)}{\Amap(-)}{}
\def\Rr(#1){R_{#1}}                        % restriction of E
\gge{Rr(-)}{\Rr(-)}{restriction of E}
\def\Cr(#1){C_{#1}}                        % restriction of E to N
\gge{Cr(-)}{\Cr(-)}{restriction of E to N}
\newcommandmacro{\Tm}{{\cal T}}{Tm}{Transfer Matrix}
\def\On(#1){{\cal I}_{#1}}
\gge{On(-)}{\On(-)}{}
\newcommandmacro{\UU}{\underline{\sqcup}}{UU}{}  
\newcommandmacro{\UUU}{\sqcup}{UUU}{}  
\newcommandmacro{\Vq}{V_Q^{\otimes n}}{Vq}{Potts config. space}
\def\bs(#1,#2){\mbox{{\Large $\ast$}}^{#1}_{#2}}  % general plumbing multiplier
\gge{bs(-,-)}{\bs(-,-)}{general plumbing multiplier}
\newcommand{\ignore}[1]{}
\gge{ignore\{-\}}{\ignore{-}}{ignore argument!}
\def\choo(#1,#2){ \left( \begin{array}{c} #1 \\ #2 \end{array} \right) } %choose
\gge{choo(-1,-)}{\choo(-1,-)}{choose}
\newcommand{\Qed}{$\Box$}%{Qed}{QED}
\gge{Qed}{\mbox{\Qed}}{QED}
\def\staq(#1){\stackrel{#1}{=}}            % stack =
\gge{staq(-)}{\staq(-)}{}
\def\stam(#1){\stackrel{#1}{\rightarrow}}  % stack ->
\gge{stam(-)}{\stam(-)}{}
\def\mat{ \left( \begin{array} }    
\def\tam{ \end{array}  \right) }
\gge{mat/tam}{...}{matrix delimiters}
\newcommand{\beq}{\begin{equation} }
\def\eql(#1){ \begin{equation} \label{#1} 
}
\newcommand{\eq}{\end{equation} }
\def\eqal(#1){\begin{eqnarray} \label{#1} }
\def\eqa{\end{eqnarray} }
\def\lab(#1){\label{#1}
}
\def\prl(#1){ \begin{pr} \label{#1} 
}

\def\leml(#1){ \begin{lem} \label{#1} 
}

\def\del(#1){ \begin{de} \label{#1} 
}
\gge{smeq\{-\}}{...}{small equation}
\gge{fneq\{-\}}{...}{very small equation}
\noglossaryignore{\newline HECKE/BLOB\newline}
\newcommandmacro{\Hnq}{H_n(q)}{Hnq}{ * freestanding symbol}
\newcommandmacro{\Hn}{H_n}{Hn}{      *-mod etc.}
\newcommandmacro{\A}{{\cal A}}{A}{}
\newcommandmacro{\Cwts}{C}{Cwts}{}
\newcommandmacro{\CA}{{\cal A}}{CA}{}

\newcommandmacro{\calA}{{\cal A}}{calA}{}
\newcommandmacro{\modi}{\mbox{Mod} }{modi}{was mod not modi!}
\newcommandmacro{\Wgen}{{\Bbb S}}{Wgen}{}
\def\ol(#1){\overline{#1}}
\newcommandmacro{\st}{\mbox{St}}{st}{}
\def\CMult(#1,#2){(#1:#2)}
\def\CM(#1,#2){( #1 : #2 )}
\def\FMult#1,#2{(#1:#2)}
\def\CF#1,#2{(#1:#2)}

\newcommandmacro{\Top}{\mbox{Top}}{Top}{}
\newcommandmacro{\Soc}{\mbox{Soc}}{Soc}{}
\newcommandmacro{\Head}{\mbox{Head}}{Head}{}
\newcommandmacro{\Filt}{{\cal F}}{Filt}{}
\newcommandmacro{\Mod}{\mbox{mod}}{Mod}{}
\newcommandmacro{\Resi}{\mbox{Res }}{Resi}{was without i!}
\newcommandmacro{\Indi}{\mbox{Ind }}{Indi}{was without i!}
\def\RR(#1,#2){R^{#1}_{#2}}   %projection after restriction (after projection)
\def\TT(#1,#2){T^{#1}_{#2}}   %translation

\newcommandmacro{\Ann}{\mbox{Ann}}{Ann}{}
\newcommandmacro{\Cen}{\mbox{Cen}}{Cen}{}
\newcommandmacro{\End}{\mbox{End}}{End}{}
\newcommandabbreviation{\semisimple}{semisimple}{semisimple}
\newcommandabbreviation{\Bratteli}{Bratteli}{Bratteli}
\newcommandabbreviation{\JBC}{Jones Basic Construction}{JBC}
\newcommandabbreviation{\pa}{partition algebra}{pa}
\newcommandabbreviation{\TM}{transfer matrix}{TM}
\newcommandabbreviation{\PM}{Potts model}{PM}
\newcommandabbreviation{\QSC}{quantum spin chain}{QSC}
\newcommandabbreviation{\Hamiltonian}{Hamiltonian}{Hamiltonian}
\newcommandabbreviation{\YS}{Young symmetrizer}{YS}
\begin{document}
\title{Tiling bijections between paths and Brauer diagrams}
\date{28 June 2010}
\author{Robert J Marsh and Paul Martin}
\thanks{This work was supported by the Engineering and Physical Sciences
Research Council [grant number EP/G007497/1]}.

\begin{abstract}
There is a natural bijection between Dyck paths and basis diagrams of the
Temperley-Lieb algebra defined via tiling.
\emph{Overhang} paths are certain generalisations of Dyck paths
allowing more general steps but restricted to a rectangle in the
two-dimensional integer lattice.
We show that there is a natural bijection, extending the above
tiling construction,
between overhang paths and basis diagrams of the Brauer algebra.
\end{abstract}
\maketitle

\vskip 0.5cm

\newcommand{\aar}{\ar@{-}}          % arrowhead
\newcommand{\headroom}{80}        % generic object height
\newcommand{\raised}{-\headroom}  % generic object offset
\newcommand{\pp}[2]{\begin{picture}(#1,\headroom)(0,\raised)
    \put(-20,0){#2} \end{picture}}
\newcommand{\ppp}[4]{\begin{picture}(#1,#2)(0,-#3) \put(-20,0){#4}
    \end{picture}}        %#1=x size,  #2=y size,  #3=y offset #4=picture content

\newcommand{\putch}{\put(-.3,-26)}
\newcommand{\poutch}{\put(-.3,-20)}
%%%%%%%%%%%%%%%%%
\newcommand{\overhang}{overhang}
\newcommand{\Dyke}{Dyck}
\newcommand{\Dyck}{Dyck}
\newcommand{\TLproj}{\Pi} %%hook for map from Brauer to TL diagrams.
\newcommand{\TLinv}{\Phi^{TL}}

% We later define \Y for the set of overhang paths and \YTL for the set
% of Dyck paths
%
%}}}

\newcounter{parano}[section]%% numbered paragraph
\renewcommand{\theparano}{\thesection.\arabic{parano}}
\newcommand{\p}{\refstepcounter{parano}\noindent\theparano.\ }

\newcounter{mysubsection}[section]%% numbered paragraph
\renewcommand{\themysubsection}{\thesection\Alph{mysubsection}}
\newcommand{\mysubsection}[1]{
\begin{center}
\refstepcounter{mysubsection}\noindent\themysubsection. \textbf{#1}
\end{center}}

\newtheorem*{theo*}{Theorem}
\newtheorem*{cor*}{Corollary} 
\newtheorem*{pr*}{Proposition} 
\newtheorem*{co*}{Corollary}
\newtheorem*{rem*}{Remark} 
\newtheorem*{lem*}{Lemma} 
\newtheorem*{claim*}{Claim}
\newtheorem*{ex*}{Example}

\newcommand{\defn}{\textbf{Definition.\ }}
\newcommand{\defcom}[1]{\textbf{Definition (#1).} \\}

%%%%%%%%%%%%%%%%%%%%%%%%%%%%%%%%%%%%%%%%%%%%%%%%%%%%%%%%%%%%%%%%%%%%%%%%
%This is Ch-intro.tex
%%%%%%%%%%%%%%%%%%%%%%%%%%%%%%%%%%%%%%%%%%%%%%%%%%%%%%%%%%%%%%%%%%%%%%%%
\renewcommand{\ul}[1]{\underline{#1}}

\section{Introduction}
%{{{ 1

Consider the double factorial sequence, given by
$S_n=(2n-1)!!=(2n-1)(2n-3)\cdots 1$.
The sequence begins:
\[
1,3,15,105,945, \ldots
\]
There are many important sequences of sets whose terms have cardinalities
given by this sequence (see, for example, entry A001147 of~\cite{sloane09}).
The `abstract' challenge is, given a
pair of such sequences, to find bijections between the $n$th terms in
each sequence that are \emph{natural} in the sense that they
can be described for all $n$ simultaneously.
We consider here Brauer diagrams (pair partitions of $2n$
objects) and \overhang\ paths (certain walks on a rectangular grid).

A striking example of a natural bijection, for the sequence of Catalan numbers,
is the bijection between Temperley-Lieb diagrams
(non-crossing pair partitions) and Dyck paths
(see e.g.~\cite{stantonwhite86}), given by `tiling'. Recall that a Dyck path
is a non-collapsing path in the upper half-plane starting at the origin in
which each step increases the $x$-coordinate by $1$ and changes the
$y$-coordinate by $\pm 1$, here with a specified end-point on the $x$-axis.
Here is an example of a tiling of a Dyck path giving rise to a
Temperley-Lieb diagram:

\begin{center}
\includegraphics[height=1in]{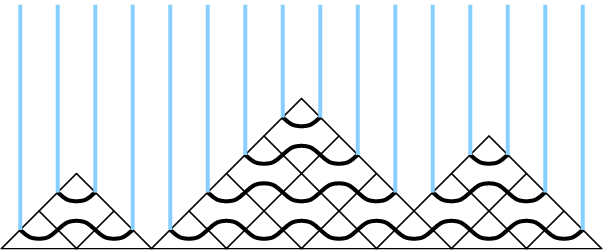}
\end{center}

See Sections~\ref{s:tilemap} and~\ref{s:TLtile} for more details.

The Dyck path basis of standard modules over the Temperley-Lieb
algebra~\cite{temperleylieb71} lends itself to the construction of
Young's orthogonal form for such modules.
The Young tableau realisation of Specht modules
plays a similar role for the symmetric group algebra and the Hecke
algebra. From this one is able to read off the `unitarisable' part of the
representation theory of the algebra in question for $q$ a root of unity
--- that is, the simple modules appearing in Potts tensor
space~\cite[\S8.2]{Martin91}. This is much harder to do using the
Temperley-Lieb diagrams themselves, where the necessary combinatorial
information is completely obscure. In fact, the Temperley-Lieb diagrams
define instead the fundamental {\em integral} form of the corresponding
modules. Therefore, the bijection between Temperley-Lieb diagrams and
Dyck paths provides a good example of an interesting bijection from
a representation theory perspective.

%}}}
%{{{ 1.0

Much progress has been made recently (see e.g.~\cite{cdm09} and
references therein) on the representation theory of
the Brauer algebra~\cite{brauer37} but an analogue of the orthogonal
form/simple module construction cited above (and described in
Section~\ref{s:TLparadigm} in greater detail) is not known. For this
reason, as a first step towards this, it is of interest to construct
a parallel bijection between overhang paths and Brauer diagrams. We do
this here.

An \emph{overhang path} is defined in the same way as a Dyck path,
except that steps in which the $x$-coordinate is decreased by $1$ and
the $y$-coordinate is increased by $1$ are also allowed.
In addition, the path is not allowed to cross the $y$-axis.
The proof that the map we construct is a bijection is nontrivial but
a flavour can be given by the following, in which a tiling of an
overhang path gives rise to a Brauer diagram:

\begin{center}
\includegraphics[height=1.3in]{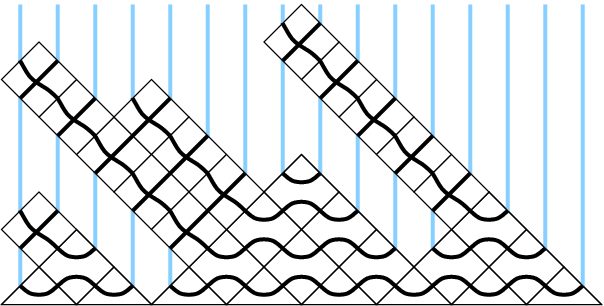}
\end{center}

See Section~\ref{s:tilemap} for the definition of the tiling map,
and sections~\ref{s:BrauertoDyck} to~\ref{s:mainresult} for the proof that
it is a bijection.

The eventual aim is to push this result on into representation theory,
as in the Temperley-Lieb case, but we
restrict here to reporting on the initial combinatorial work necessary.

The non-crossing pair partitions (Temperley-Lieb diagrams) are a subset of 
the set of general pair partitions. Dyck paths are a subset of the set
of \overhang\ paths. With this in mind we require that our bijection agrees
with the Temperley-Lieb/Dyck path correspondence when restricted to
Temperley-Lieb diagrams.

There is in fact another bijection between Brauer diagrams and \overhang\ 
paths that is relatively easy to construct, but it does not preserve the
Temperley-Lieb/Dyck path correspondence in the above sense.
We describe this simpler correspondence in Section~\ref{s:easymap}.

The article is organised as follows.
In Section~\ref{s:defnsetc},
we discuss Dyck paths and \overhang\ paths and their properties.
In Section~\ref{s:brauer},
we recall Brauer diagrams and define some simple notions
on such diagrams which will be useful later.
In Section~\ref{s:tilemap} we define a
\emph{tiling} map from \overhang\ paths to Brauer diagrams.
In Section~\ref{s:TLtile} we recall a tiling-type bijection between
Dyck paths and Temperley-Lieb diagrams.
In sections~\ref{s:BrauertoDyck} to~\ref{s:mainresult}
we show that the map in Section~\ref{s:tilemap} has an inverse,
thus proving our main result, Theorem~\ref{th:main}, that there is a bijection
between \overhang\ paths and Brauer diagrams which extends the bijection
described in Section~\ref{s:TLtile}. In Section~\ref{s:example}
we give an example.
In Section~\ref{s:easymap} we describe the simpler bijection between
Brauer diagrams and \overhang\ paths (which does not extend the
tiling map in the Temperley-Lieb case).
Finally, we explain some of our motivation in
terms of the orthogonal form construction in the Temperley-Lieb/Dyck
path setting in Section~\ref{s:orthogonalform}.

We would like to thank M.\ Grime for bringing to our
attention a certain notion of paths in the plane
(we refer to them here as \emph{\overhang\ paths of degree $n$};
see~\ref{p:overhangdef}), and also for his initial question which motivated
us to start work on this article. He mentioned to us that it was known
that the number of overhang paths of degree $n$ coincides with the number
of Brauer diagrams of degree $n$ (for formal reasons: the generating
functions are identical) and asked the question as to whether this could
be proved concretely.

After we completed work on this article, we learnt of the
article~\cite{rubey}, which also gives a bijection between overhang paths and
Brauer diagrams. This bijection is different from both of the bijections
we define here, and we do not know a way of defining it using tilings.

We remark that there are a number of other examples of sets in natural
bijection with Brauer diagrams. As well as those in entry A001147
of~\cite{sloane09}, there are examples in~\cite{dalemoon93}.
For information on bijections between Brauer diagrams
(and more general partitions) and tableaux and pairs of walks,
we refer to~\cite{cddsy07,halversonlewandowski05,martinrollet98,marshmartin,sundaram86,sundaram90,terada01}.
We also remark that the article~\cite{bakerforrester01} gives a bijection
between fixed-point free involutions of a set of size $2n$
and certain sets of tuples of non-intersecting walks on the natural
numbers arising in statistical mechanics (the random-turns model of vicious
random walkers).

%}}}
%}}}

%%%%%%%%%%%%%%%%%%%%%%%%%%%%%%%%%%%%%%%%%%%%%%%%%%%%%%%%%%%%%%%%%%%%%%%%%%%%%
%End of Ch-intro.tex
%%%%%%%%%%%%%%%%%%%%%%%%%%%%%%%%%%%%%%%%%%%%%%%%%%%%%%%%%%%%%%%%%%%%%%%%%%%%%

\section{\overhang\ paths}
\label{s:defnsetc}
%{{{ defns etc

\p
Consider the semi-infinite rectangle $R\subseteq \mathbb{R}^2$
with base given by the line segment from $(0,0)$ to $(n,0)$ and sides
$x=0$ and $x=n$.
Let $R_{\mathbb{Z}}$ denote the set of integral points $(a,b)$ in this
rectangle.
We consider steps between points in $R_{\mathbb{Z}}$ of the following form:

\begin{description}
\item[($1$)] $(a,b) \rightarrow (a+1,b+ 1)$, or 
\item[($2$)] $(a,b) \rightarrow (a+1,b- 1)$, or
\item[($2'$)] $(a,b)\rightarrow (a-1,b+1)$.
\end{description}

\newcommand{\Y}{\mathcal{G}}  %% hook for grime paths
\newcommand{\YTL}{\mathcal{G}^{TL}} %%hook for Dyck paths

\p 
We define a \emph{\Dyck\ step} to be a straight line path of form $(1)$ or
$(2)$, and a \emph{\overhang\ step} to be a straight line path of form $(2')$.

\begin{figure}
\[
\includegraphics[width=2in]{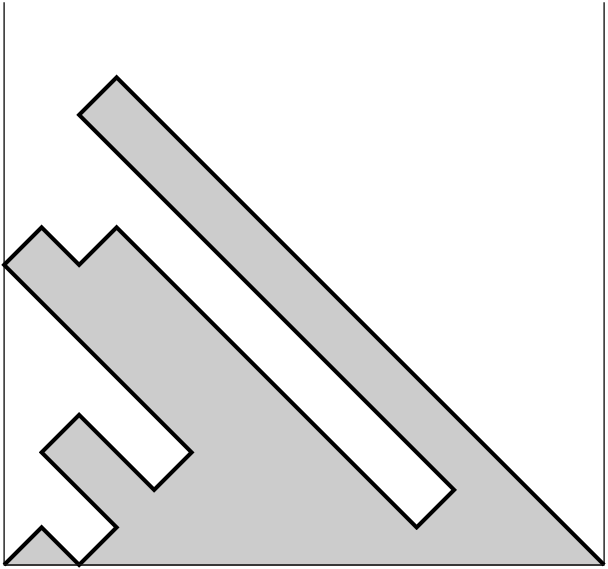}
\]
\caption{An \overhang\ path from $(0,0)$ to $(8,0)$.}
\label{fig:examplegrime}
\end{figure}

\p \label{p:overhangdef}
A \emph{path} in $R_{\mathbb{Z}}$ is a sequence of steps between vertices of
$R$. It is said to be \emph{noncollapsing} if it does not visit any vertex
more than once. In particular, a \emph{Dyck path}
(respectively, \emph{\overhang\ path})
is a noncollapsing path starting at $(0,0)$ and consisting of Dyck
(respectively, Dyck or \overhang) steps.
We shall restrict our attention
to Dyck or \overhang\ paths which end at $(2n,0)$ for some $n\in\mathbb{N}$;
such paths will be said to have \emph{degree} $n$.
Let $\YTL_n$ (respectively, $\Y_n$) denote the set of all
\Dyck\ (respectively, \overhang) paths of degree $n$.
For an example of an \overhang\ path of degree $8$, see
Figure~\ref{fig:examplegrime} (the shading in the figure will be explained
in~\ref{p:lowerregion}).

%}}}
%{{{ fig 

\p
There is an injective map from paths to finite sequences of elements from the set
$\{1,2,2'\}$ given by writing a path as its sequence of steps. For example,
$$\Y_2=\{1122,12'1222,1212\}.$$

%}}}
%{{{ remarks

\p \label{p:lowerregion}
A path $p \in \Y_n$, together with the $x$-axis with
the interval between $(0,0)$ and $(2n,0)$ removed, partitions the plane
into two regions. The intersection of these regions with $R$ will be
referred to as the \emph{upper region} and the \emph{lower region}
of $p$ respectively. (In the example in Figure~\ref{fig:examplegrime},
the lower region is shaded.)

\p
We define a partial order on $\Y_n$ by setting $p<q$ if 
the lower region of $p$ is contained in the lower region of $q$.
Thus, the lowest path is
$$p_0=121212...12.$$

\p
If $p<q$, we shall write $q/p$ for the `skew' diagram --- the lower region of
$q$ not in the lower region of $p$.

\begin{figure}
\[
\includegraphics[width=2.5in]{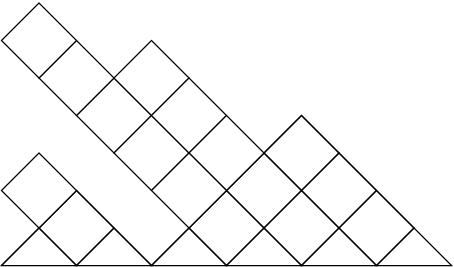}
\]
\caption{\label{fig:tilingexample} Tiling an \overhang\ path.}
\end{figure}

\p
We will consider the lower region of $p$ not in the lower region
of $p_0$ to be tiled with diamond tiles, and we will consider the lower
region of $p$ intersecting the lower region of $p_0$ to be tiled with
half-diamond tiles. For an example, see Figure~\ref{fig:tilingexample}.

%}}}
%{{{ degree etc

\begin{lem*} \label{lem:grimecardinality} 
Let $n\in\mathbb{N}$. Then $|\Y_n|=(2n-1)!!$.
\end{lem*}

{\bf Proof:}
Given a sequence $\mathbf{r}=(r_0,r_1,\ldots ,r_{n-1})$ of integers satisfying
$0\leq r_k\leq 2k$ for $0\leq k\leq n-1$, we can form an \overhang\ path
in the following way. Start with the path $p_0$ described above.
Then, for each $k$, add a rectangle $R_k$ to the $p_0$ with
vertices $A_k=(2k,0),B_k=(2k+1,1),C_k=(2k-r_k+1,r_k+1)$ and $D_k=(2k-r_k,r_k)$
(this can be considered as a pile of $r_k$ diamonds piled up to the left of
the step from $(2k,0)$ to $(2k+1,1)$ of $p_0$).
The upper boundary of the union of these rectangles consists of steps of form
(1) (corresponding to a line segment $D_kC_k$), form (2)
(corresponding to part of a line segment $C_kB_k$ in the case where
$r_k\geq r_{k+1}$ or $k=n$), or form (2') (corresponding to part of a line
segment $A_kD_k$ in the case where $r_k\leq r_{k-1}$).
Hence this forms the lower region of an \overhang\ path.

\begin{figure}
\[
\includegraphics[width=2in]{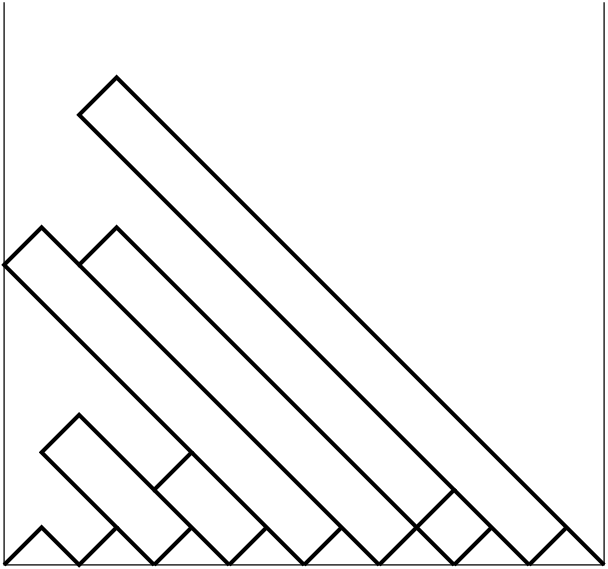}
\]
\caption{Constructing an \overhang\ path from rectangles.}
\label{fig:rectangularstrips}
\end{figure}

Conversely, given an overhang path, steps of form $(2)$ or $(2')$ from
$(a,b)$ do not change the sum $a+b$, while a step of form $(1)$ increases it
by $2$. It follows that any \overhang\  path must contain precisely $n$ steps
of form $(1)$. By considering the diamond tiles down and to the right of these
steps, we see that the path must be of the above form. It is clear that we
now have a one-to-one correspondence between overhang paths and tuples
of integers as above. The result follows. For an example, with
$\mathbf{r}=(0,0,3,2,8,8,1,12)$, see Figure~\ref{fig:rectangularstrips}.~$\Box$

\p \label{p:root}
For each $p \in \Y_n$ there is a unique maximal path $t\leq p$
that only uses \Dyck\ steps. We call this the \emph{root} Dyck path
(or just the root) of $p$. For example, in the introduction,
the Dyck path example is the root of the overhang path example.

\p
For $p \in \Y_n$ and $q \in \Y_m$ the side-by-side concatenation
$p*q$ of $p$ and $q$ is a path in $\Y_{n+m}$:
\[
* : \Y_n \times \Y_m \rightarrow \Y_{n+m}.
\] 
Note that not every path in $\Y_{n+m}$ that passes through $(2n,0)$
arises in this way.

\p
An element of $\Y_n$ is said to be \emph{prime} if it cannot be expressed
non-trivially in the form $a*b$.

%}}}

\section{Brauer diagrams}
\label{s:brauer}
%{{{ brauer
%{{{ sec 2

\newcommand{\un}[1]{\underline{#1}}

\begin{figure}
\[
\includegraphics[width=4in]{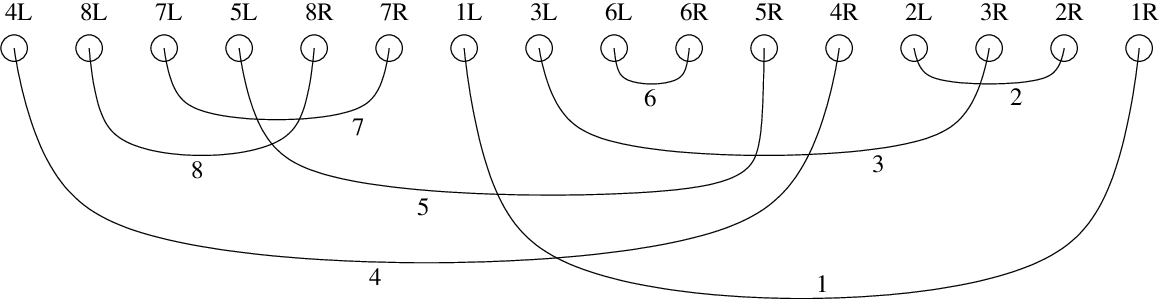}
\]
\caption{\label{fig:brauerexample} A Brauer diagram with arc labels.}
\end{figure}

\p
Given a finite set $S$, a \emph{pair partition} of $S$ is a partition
of $S$ into subsets of cardinality $2$.
A \emph{Brauer diagram} of degree $n$ is a picture of a pair partition of
$2n$ distinct vertices arranged on the boundary of the lower
half-plane. The two vertices in each part of the pair partition are joined
by an arc in the lower half-plane.
Two Brauer diagrams are identified if their underlying vertex pair partitions
are the same. Let $J_n$ denote the set of all Brauer diagrams of degree $n$.
See Figure~\ref{fig:brauerexample} for an example.
The additional arc and vertex labels are explained below. 

\p We remark that Brauer diagrams are often defined using $2n$ vertices on
the boundary of a disk or in a horizontal rectangle, with $n$ vertices along
the top and $n$ vertices along the bottom, but we shall not consider such
representations here.

\p
By a \emph{partial Brauer diagram}, we mean a Brauer diagram, but with
the extra possibility that parts of cardinality $1$ are also allowed.
We denote by $J_n^l$ the set of partial Brauer diagrams containing $n$
pairs and $l$ singletons (and thus a total of $2n+l$ vertices).

\p
See Figure~\ref{fig:rightagree}(a) for an example of a partial Brauer diagram
which is not a Brauer diagram.
We remark that a partial Brauer diagram can be completed on the left by
adding another partial Brauer diagram to the left with the same number of
singletons, and then pairing up the singletons in the first diagram with
those in the second. Note that such a completion is in general not unique.

\p
A \emph{TL diagram} (or Temperley-Lieb diagram) is a Brauer diagram without
crossings. We shall write $J_n^{TL}$ for the subset of $J_n$ consisting of
TL diagrams.

%}}}
\newcommand{\LH}{left-hand}
\newcommand{\RH}{right-hand}
%{{{ arc labels (incl fig)

\p 
\defcom{(Right-)standard arc labelling}
Let $D$ be a partial Brauer diagram. We number the vertices of
$D$ which are right-hand ends of arcs or singletons, in order from right
to left. A vertex $k$ which is the right-hand end of an arc gets labelled
$kR$, and we label the other end of the arc $kL$. Sometimes we will label
the arc with endpoints $kL$ and $kR$ with the number $k$.

For an example, see Figure~\ref{fig:brauerexample}.

\p
We define similarly a {\em left-standard} labelling, which again
numbers from right to left, but according to the order of the {\em left}-hand
endpoints of arcs (and singletons as before).

\p
Later we will use the pair $(a(i),i)$ of left and right-standard labels
for an arc in a fixed diagram $D$.
That is, if $i$ is the right-standard label of an arc, then $a(i)$ will be 
the left-standard label of the same arc. 

\p
We will not need the, perhaps more natural, orderings from left to
right. This handedness comes from the handedness of the \overhang\ diagrams
that we chose.

%}}}
%{{{ subdiagram

\p
To each arc $i$ (in the right-standard labelling) of a diagram $D$ we may
associate an {\em arc (left) subdiagram} $D^i$ of $D$.
This is the collection of arcs whose \RH\ vertex is
strictly contained within arc $i$ (i.e.\ the interval from $iL$ to $iR$),
together with their endpoints. We retain the initial (right-standard)
labelling of the vertices inherited from $D$.

\begin{figure}
\[
\includegraphics[width=3.4in]{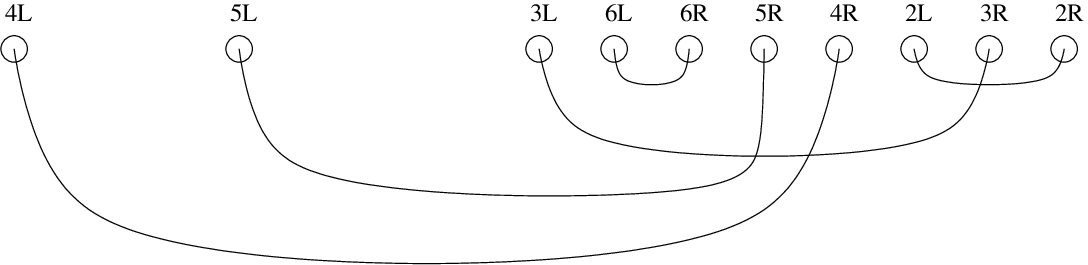}
\]
\caption{\label{fig:arcsubdiagramexample} The arc subdiagram $D^1$.}
\end{figure}

\p
\textbf{Example:} Let $D$ be the diagram in Figure~\ref{fig:brauerexample}
above.
The arc subdiagram $D^1$ is shown in Figure~\ref{fig:arcsubdiagramexample}.

%}}}
%{{{ monoidal stuff

\def\ne(#1){[ #1 ]} %% the nest bracket notation

\p
To any diagram $D$ in $J_n$ we may associate a diagram in $J_{n+1}$,
denoted $\ne(D)$, which is the diagram obtained from $D$ by adding a new
vertex at each end of $D$, and an arc between them. 

Similarly for $D \in J_n$ and $D' \in J_m$ we will understand by 
$DD' \in J_{n+m}$ the diagram obtained by simple side-by-side
concatenation. 

\p
We shall call a diagram {\em prime} if it cannot be expressed non-trivially in
the form $D=D_1 D_2$. 
(This is a different definition of prime than has been used
  elsewhere, e.g.~\cite{lickorish97,martinsaleur94}).
Note that if a diagram $D$ is Temperley-Lieb and prime then it can be
expressed in the form $D = \ne(D')$. 

%}}}
%}}}

\section{The tile map}
\label{s:tilemap}
%{{{ tiles

\p
There is a map from \overhang\ paths to Brauer diagrams
\[
\Psi : \Y_n \rightarrow J_n
\]
defined by replacing each `blank' tile with a patterned tile.
Tiles in the root Dyck path of $p\in \Y_n$ are replaced using the
following rules: 
\[
\includegraphics[width=1.5in]{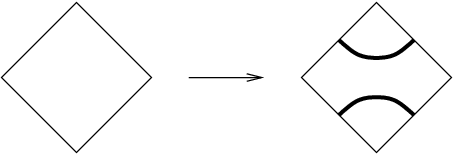}
\]
or
\[
\includegraphics[width=1.5in]{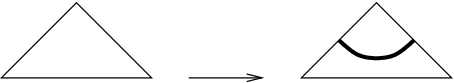}
\]
 Tiles in the lower region of $p$ but above the root of $p$ are replaced using
the following rule:
\[
\includegraphics[width=1.5in]{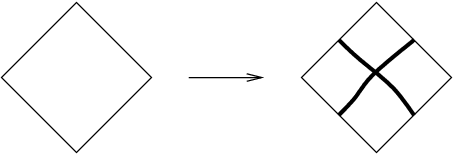}
\]

A horizontal line above the overhang path is fixed (the ``top'' of the
diagram). Strands are then connected together with vertical segments joining two
ends, or joining an end with the top of the diagram. This can also be
realised by continuing the tiling into the upper region of the path
(up to the horizontal line), using half-tiles on the boundary, and using the
following tiling rules for the new tiles:

\[
\includegraphics[width=1.5in]{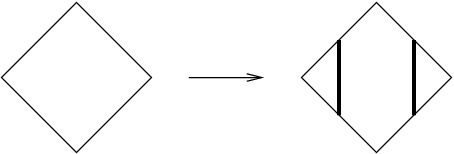}
\]

\[
\includegraphics[width=1in]{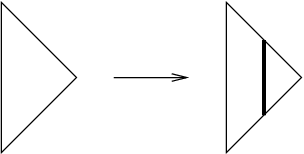}
\]

\[
\includegraphics[width=1in]{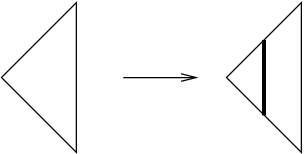}
\]

\[
\includegraphics[width=1.5in]{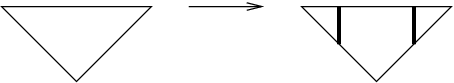}
\]

\begin{figure}
\[
\includegraphics[width=4in]{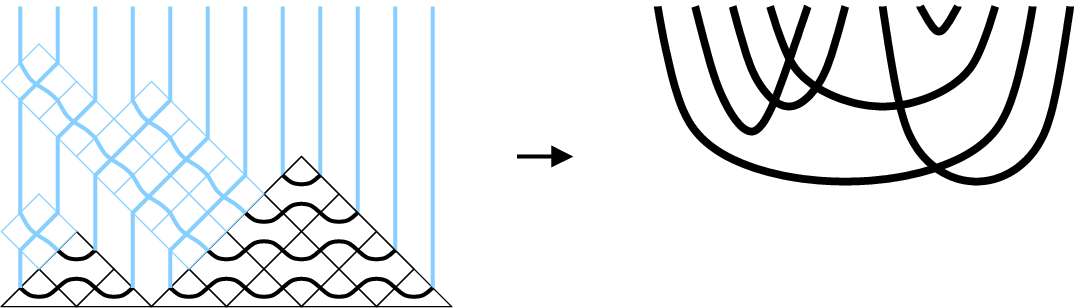}
\]
\caption{Computation of the Brauer diagram $\Psi(p)$ for the path $p$ in
Figure~\ref{fig:tilingexample}.}
\label{fig:mapexample}
\end{figure}

For an example (with the tiles in the upper region omitted for clarity),
see Figure~\ref{fig:mapexample}.

\p
We note that the patterned tiling of the lower part of $p$ in the above
construction can be regarded as a \emph{pipe dream}~\cite{fominkirillov96}
(also known as an \emph{rc-graph}~\cite{bergeronbilley93}).
In general it will be non-reduced, i.e.\ two arcs may cross
twice in the resulting configuration (see the introduction for
an example of this).

\p
By the construction of $\Psi$, we have:

\begin{lem*} \label{lem:coms}
The map $\Psi$ commutes with side-by-side concatenation: 
\[
\Psi(a * b) = \Psi(a) \Psi(b) 
\]
for all $a,b\in \Y_n$.
\Qed
\end{lem*}

%}}}
\section{The Temperley-Lieb case}
\label{s:TLtile}
%{{{ TL

\p
Note that the map $\Psi$ has image within the set of TL diagrams when 
restricted to the set of roots, given by `tiling':
\[
\Psi|_{\YTL_n}: \YTL_n \rightarrow J_n^{TL}
\]
See e.g.~\cite{stantonwhite86}.

\p
The inverse of the restricted map is also well known. A convenient
in-line representation of a TL diagram $D$ is to read from left to right
and to replace each vertex that is the left hand end of an arc with an
open bracket, "(", and to replace each vertex that is the right hand end of
an arc with a close bracket, ")". 
It is clear that this gives rise to a well-nested sequence of
brackets. Replacing each "(" with a $1$ and each ")" with a $2$ we
obtain the in-line sequence for a \Dyck\ path, call it  $\TLinv(D)$. 

\p
By construction, we have the following:

\begin{lem*} \label{lem:coms2}
The map $\TLinv$ commutes with side-by-side concatenation:
\[
\TLinv(D  D') = \TLinv(D) * \TLinv(D') 
\]
for all TL diagrams $D,D'$.
\Qed
\end{lem*}

%}}}
%{{{ more TL

\begin{lem*} \label{lem:sx}
The map $\TLinv$ is the inverse of $\Psi|_{\YTL_n}$.
\end{lem*}

\noindent {\em Proof:}
This is implicit in~\cite{abf84} (see~\cite{Martin91}),
but we include a proof for the convenience of the reader.
We show that for all TL-diagrams $D$, $\Psi(\TLinv(D))= D$. 
We do this by induction on $n$, with $n=0$ as base. 
Suppose that the result is true for smaller $n$.
If $D$ is a TL diagram of degree $n$, suppose first that $D$ has an arc
joining vertices $1$ and $2n$. Let $D'$ be the TL-diagram obtained by
removing this arc. By induction, $\Psi(\TLinv(D'))=D'$. It follows that
$\Psi(\TLinv(D))=D$, since $\TLinv(D)$ is the same as the Dyck path $\TLinv(D')$
except that an extra step $1$ at the start and an extra step $2$ at the
end have been added.

If $D$ has no arc joining $1$ and $2n$ then it is of the form $D_1D_2$
where $D_1$ and $D_2$ are non-empty TL diagrams. By the inductive
hypothesis, $\TLinv(\Psi(D_i))=D_i$ for $i=1,2$, and it follows that
$\TLinv(\Psi(D))=D$. The result follows by induction.

It is well known that the cardinalities of $\YTL_n$ and $J_n^{TL}$
are the same (given by the $n$th Catalan number), so the result follows.
\Qed

%}}}

%{{{ entile
\section{A map from Brauer diagrams to Dyck paths}
\label{s:BrauertoDyck}
%{{{ intro
\p
Our ultimate aim is to define a map
\[
\Phi : J_n \rightarrow \Y_n
\]
and show that it is inverse to $\Psi$. The difficulty is that the overhang
path corresponding to a Brauer diagram may be hard to find. In the
example given in the introduction, the realisation of the Brauer diagram
obtained from the overhang path is not the simplest one: it contains more
crossings than are necessary (one of the strings crosses one of the others
twice: both crossings could be removed). Our approach will be to find
first what will turn out to be the root of the desired overhang path and then
add extra tiles to the corresponding Dyck path in order to give the required
crossings.

Thus we will first of all define a map
$\TLproj:J_n \rightarrow J^{TL}_n$ associating a Temperley-Lieb diagram
to each Brauer diagram. In the next section, we shall see that this gives us
a useful labelling of each Brauer diagram. We will then study the properties
of this labelling. This will give us control of the crossings and allow us to
define $\Phi$. (The definition of $\Phi$ appears in Definition~\ref{def:inversemap} and the main result is Theorem~\ref{th:main}.)
%}}}
%{{{ chain

\p
The {\em (right) chain} $ch(D)$ of arcs of $D \in \Y_n$ is the sequence
$a_1,a_2,\ldots $ of arc labels of $D$ such that $a_1=1$ and $a_i$ 
(if it exists) is the arc label of the first \RH\ end vertex to occur
moving from right to left from the \LH\ end vertex of the arc with
label $a_{i-1}$. 

\p
\textbf{Example:} for the diagram in Figure~\ref{fig:brauerexample},
we have $ch(D) = (1,7)$ (in the right standard labelling);
while $ch(D^1) = (2,4)$ (borrowing the same labelling). 

Note that the set of arcs in the chain $ch(D)$, together with the
sets of arcs in the arc subdiagrams for the arcs in the chain 
(that are clearly disjoint sets), form the complete set of arcs of $D$. 

%}}}
%{{{ tree

\p
\defcom{Right chain  %%order on arcs
tree of $D$}
Fix a diagram $D$.
Firstly, for each arc $i$ of $D$ define a planar rooted tree with root $i$ and
other vertices the chain arcs of $D^i$ arranged in right chain order, right
to left, at tree distance 1 from the root.
For example, for the arc $1$ in Figure~\ref{fig:brauerexample}, we obtain:
\[
\includegraphics[width=1.4in]{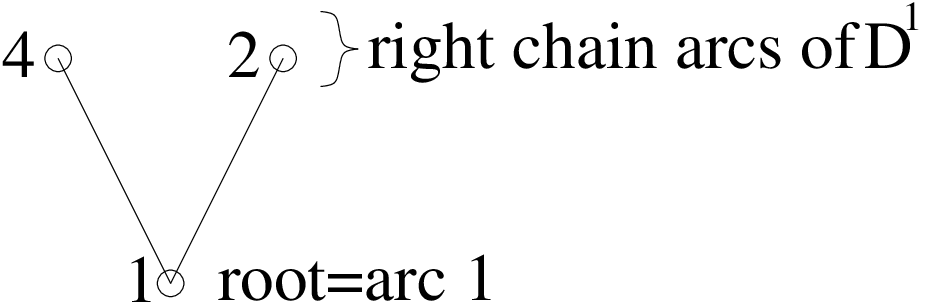}
\] 
Note that the second row contains the right chain arcs of $D^1$ given above.
Define a planar rooted tree $\tau_R(D)$ with vertices the arcs of $D$ together
with a root vertex $\emptyset$. 
The tree is obtained by concatenating the planar rooted trees for the
right chain arcs of $D$ in the obvious way, setting $D^{\emptyset}=D$ to
include the root. We call this tree the \emph{right chain tree of $D$}.
We have thus defined a map $\tau_R$ from $J_n$ to planar rooted trees. 

\begin{figure}
\[
\includegraphics[width=2.3in]{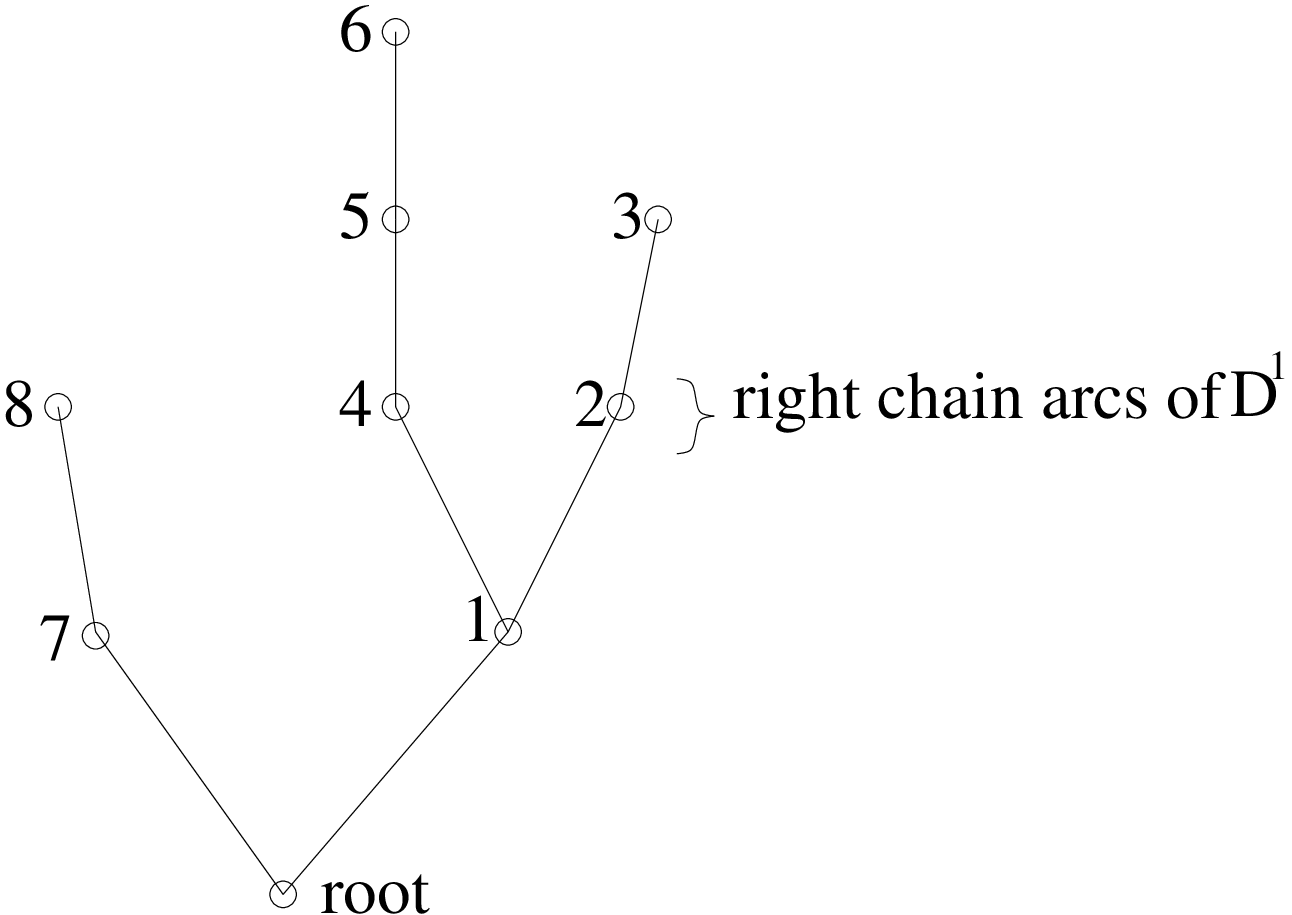}
\] 
\caption{An example of a right-chain tree.}
\label{fig:plantree12}
\end{figure}

\p
\textbf{Example}: The right-chain tree for our example $D$ above is shown in Figure~\ref{fig:plantree12}.

\begin{figure}
\[
\includegraphics[width=4in]{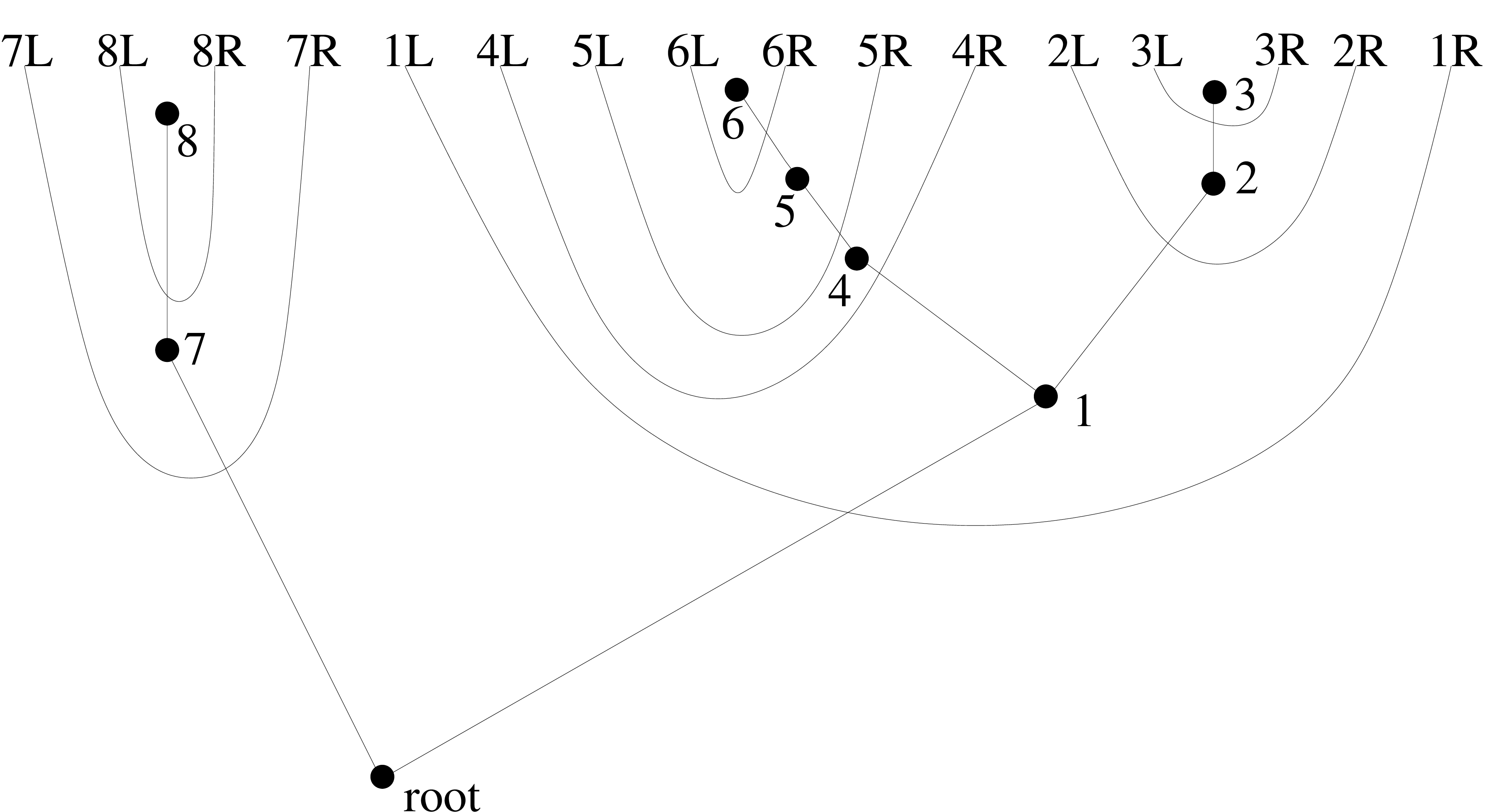}
\] 
\caption{The geometric dual of a planar rooted tree.}
\label{fig:geometricdual}
\end{figure}

\p
Let $\gamma$ denote the usual geometric dual map from planar rooted trees to
TL diagrams. Thus, for a planar tree $T$, each arc of the TL diagram
$\gamma(T)$ passes through a unique edge of $T$.
The dual TL diagram for the above example is shown in Figure~\ref{fig:geometricdual}.

Combining, we have a map
\[
\TLproj:=\gamma \circ \tau_R : J_n \rightarrow J_n^{TL}.
\]

\p
Note that the right standard labelling of the arcs of $D$
coincides with the labelling of the vertices of $\tau_R(D)$ in order of
first meeting, moving counterclockwise around the tree from the root.
An example of this can be seen in Figure~\ref{fig:geometricdual}, where
the arc ends are given their right standard labels.

\begin{figure}
\[
\TLinv(\TLproj(D)) = \includegraphics[width=3.4in]{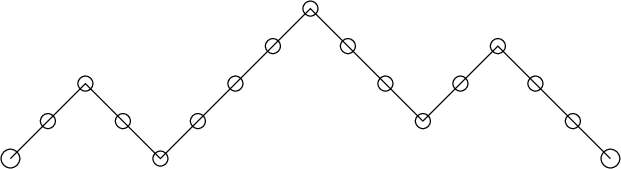}
\]
\caption{Example of a Dyck path associated to a Brauer diagram.}
\label{fig:br00105g}
\end{figure}

Note also that applying the map $\TLinv$ to $\TLproj(D)$ gives a
\Dyck\ path. So we can also associate a \Dyck\ path to each Brauer diagram.
For our example the \Dyck\ path is shown in Figure~\ref{fig:br00105g}.

\begin{figure}
\[
\includegraphics[width=4.5in]{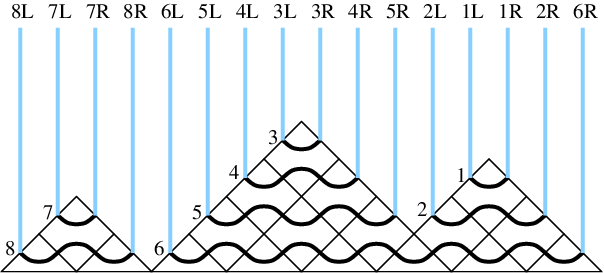}
\]
\caption{Left-standard labelling of a Temperley-Lieb diagram and
labelling of steps of form (1) in the corresponding Dyck path.}
\label{fig:br00105gTL}
\end{figure}

\p
Note that the left-standard labelling of arcs in a TL diagram induces
a labelling for steps of form (1) in the associated
\Dyck\ path, whereby each such step is given the label of the arc passing
through it. See Figure~\ref{fig:br00105gTL} for an example.

%}}}
%{{{ finish

\section{Secondary Arc Labels}
\label{s:secondary}

\p
In this section we will show how the right-standard labellings of the arcs of 
$D$ and $\Pi(D)$ can be used to obtain a new labelling (which we call secondary
labelling) of the arcs of $D$, by transferring the left-standard
labelling of $\Pi(D)$ to $D$. We shall see later how the ordering on the
arcs determined by this secondary labelling can be used to uncross the arcs
of $D$ to get $\Pi(D)$; this is a key notion in the construction of the map
$\Phi$.

\p
Fix a Brauer diagram $D$.
Each arc of the TL diagram  $\TLproj(D)$ has a pair $(a(i),i)$ of left and
right-standard labels. 
Thus for each right-standard label $i$ there is a corresponding
left-standard label $a(i)$. 

\begin{figure}
\[
\includegraphics[width=4in]{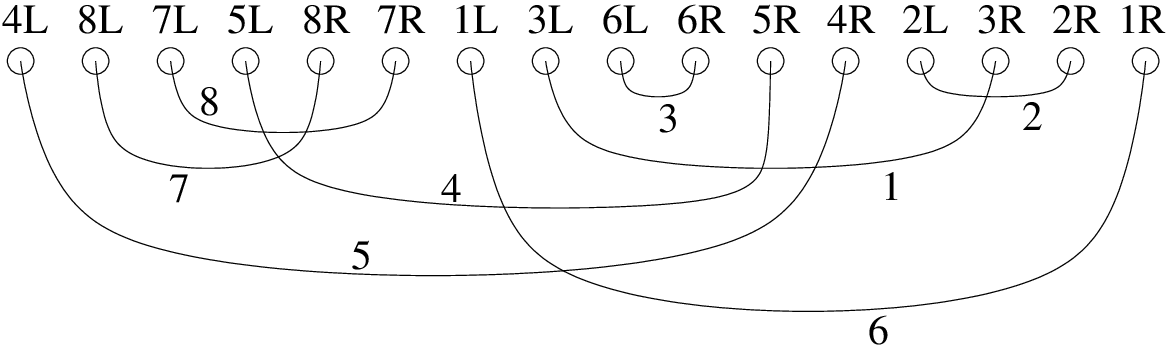}
\]
(a) A Brauer diagram $D$, with right-standard labels (at the top of the
diagram) and secondary labels (on the arcs themselves).
\vskip 0.4cm
\[
 \includegraphics[width=4in]{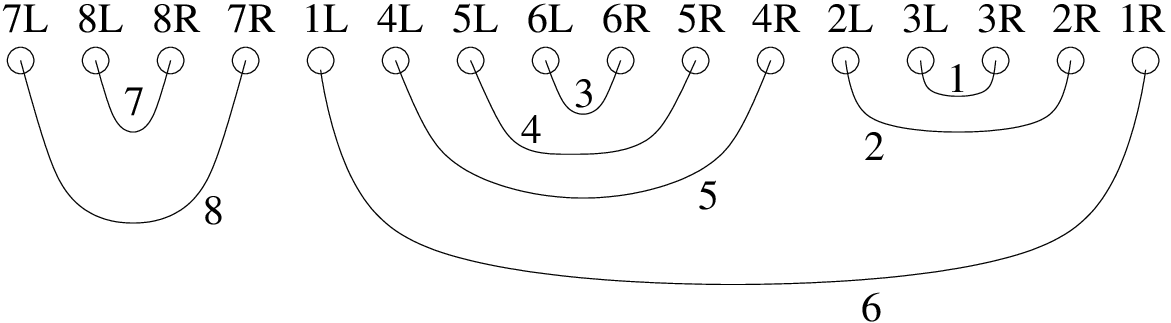}
\]
(b) The TL diagram $\TLproj(D)$ with right standard arc
 labels (on vertices), and left standard arc labels (on arcs).
\caption{A Brauer diagram $D$ and the corresponding Temperley-Lieb diagram $\TLproj(D)$, with labellings.}
\label{fig:phi'tl}
\end{figure}

For example, in the TL diagram in Fig~\ref{fig:phi'tl}(b), 
we have $a(3)=1$.

%}}}
%{{{ secondary label (relabel)

\p
\defcom{Secondary arc label}
For each arc $i$ in $D$ there is an arc in $\TLproj(D)$ 
with the same right-standard label. 
We call this association between arcs of $D$ and arcs
of $\TLproj(D)$ the `right-correspondence'. 
We now associate a new `secondary' label 
to each arc in $D$ --- the left-standard label for the right-corresponding
arc in $\TLproj(D)$. 

\p
For example, the secondary-labelling for the diagram $D$ in
Figure~\ref{fig:brauerexample} is shown in Figure~\ref{fig:phi'tl}(a).
The labels at the top of the diagram are the right-standard labels,
and each arc has been given its left-standard label.
We remark that if $D \in J_n^{TL}$ then $\TLproj(D)=D$, so its secondary and
left-standard labels coincide. 

%{{{ fig \label{fig:phi'tl}

%}}}

\p
\defcom{Right-agreement}
Let us say that two diagrams \emph{right-agree} up to a given vertex $x$
if there is a partial Brauer diagram on that vertex and the vertices to the
right of it which can be completed on the left to either of the two diagrams.
If, in addition, the two diagrams do not right-agree up to the vertex
immediately to the left of $x$, we shall say that they
\emph{maximally right-agree} up to $x$.

\begin{figure}
\[
\includegraphics[width=1.5in]{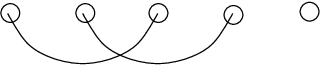}
\]
(a) A partial Brauer diagram.
\[
\includegraphics[width=3in]{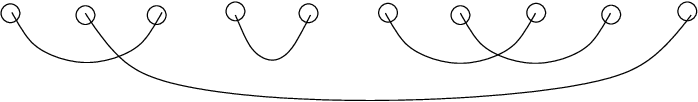}
\]
(b) A Brauer diagram.
\[
\includegraphics[width=3in]{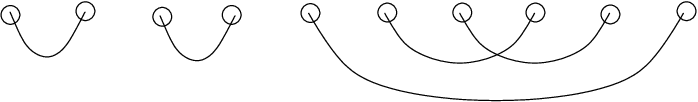}
\]
(c) A second Brauer diagram.
\caption{An example of two Brauer diagrams which right agree: each
is a completion of the the partial Brauer diagram (a).}
\label{fig:rightagree}
\end{figure}

\p
See Figure~\ref{fig:rightagree} for an example. Diagrams (b) and (c)
right-agree up to the fifth vertex from the right since the partial
Brauer diagram (a) can be completed to either of them. It is clear that
in fact the two diagrams maximally right-agree up to this vertex. 
Note also that $\TLproj(D)$ and $D$ in Figure~\ref{fig:phi'tl} maximally
right-agree up to the third vertex from the right (labelled $3R$ in both
diagrams).

%}}}
%{{{ lemmas

%\begin{lemnum}
\p \textbf{Lemma.} \label{lem:recurs1}
\emph{Suppose that $D$ and $\TLproj(D)$ right-agree up to a given vertex $x$.
Suppose that there is an arc of $D$ in the agreeing part. Then the
right-corresponding arc in $\TLproj(D)$ is also in the agreeing part
(indeed, this is the same pair of vertices in the pair partition).
These arcs have the same secondary label.
Furthermore, the set of these agreeing labels (if any) is of the form
$\{1,2,\ldots ,r\}$ for some $r$.}
%\end{lemnum}

{\em Proof:}
Let $E$ be the partial Brauer diagram which can be completed to either $D$
or $\TLproj(D)$ (on the vertex $x$ and all vertices to its right).
It follows from the definitions that the right-standard labels
on the vertices and the left-standard labels on the arcs of $E$
are the same in either completion.

Thus, on completion, an arc in $E$ gives rise to a right-corresponding pair
of arcs in $D$ and $\TLproj(D)$. By the definition of the secondary labels,
the arc in $D$ will have secondary label equal to the left-standard label
of the (right-corresponding) arc in $\TLproj(D)$. 

The consecutive property follows immediately from the definition of
left-standard labels.
\Qed

%}}}

%{{{ incomparable lemma
%\begin{lemnum}
\p \textbf{Lemma.} \label{lem:incomparable}
\emph{Let $D\in J_n$, and suppose that the left hand end of the arc with
right-standard labels $aL$ and $aR$ is to the right of the right-hand end of
the arc with right-standard labels $bL$ and $bR$.
Let $i$ be the secondary label of the former arc and $j$ the secondary
label of the latter arc.
Then in $\tau_R(D)$, $i$ and $j$ are not descendants of each other.}
%\end{lemnum}

\emph{Proof:} By the definition of $\tau_R(D)$, the descendants of
$j$ arise from (some of) the arcs whose right hand
end lies between the ends of $i$, and cannot
include $j$ by assumption. Similarly $i$ cannot
be a descendant of $j$.
\Qed

%}}}

%{{{ left hand ends lemma

\section{The relationship between $D$ and $\Pi(D)$.}
\label{s:relationship}

\p
In this section we study the relationship between $D$ and $\Pi(D)$; we shall
use these results to define $\Phi$ in the next section.

%\begin{lemnum}
\p \textbf{Lemma.} \label{lem:lefthandends}
\emph{Let $D\in J_n$, and suppose that $D$ and $\TLproj(D)$ maximally right-agree up to
vertex $x$. Then the vertex $y$ immediately to the left of $x$ is the left
hand end of an arc in both $D$ and $\TLproj(D)$.}
%\end{lemnum}

\emph{Proof:}
Case (I): If $y$ were the right-hand end of an arc in both $D$ and $\TLproj(D)$
then $D$ and $\TLproj(D)$ would right-agree up to vertex $y$, a contradiction
to the assumption. We complete the proof by ruling out the two remaining
undesirable configurations, i.e.\ the configurations in which $y$ is a
left-hand end of an arc only in $D$ or only in $\TLproj(D)$.

Case (II): Suppose that $y$ is the left-hand end of an arc in $D$ but the
right-hand end of an arc in $\TLproj(D)$.
Let $aL$ and $aR$ be the right-standard labels of the arc in $D$ incident
with $y$ in $D$ and let $bL$ and $bR$ be the right-standard labels of the arc
incident with $y$ in $\TLproj(D)$. Note that $aR$ must be $x$ or to its
right in $D$, so the vertex with right-standard label $aR$ in $\TLproj(D)$
must also be $x$ or to the right of $x$, since $D$ and $\TLproj(D)$ right-agree up
to vertex $x$, using Lemma~\ref{lem:recurs1}.
For the same reason, the vertex with right-standard label $aL$ in
$\TLproj(D)$ must be to the left of $x$, since this is so in $D$.
Since vertex $y$ is labelled $bR$ and the arcs in $\TLproj(D)$ do not
cross, the vertex with right-standard label $aL$ in $\TLproj(D)$ must be
to the left of $bL$.

\begin{figure}
\[
\includegraphics[width=3in]{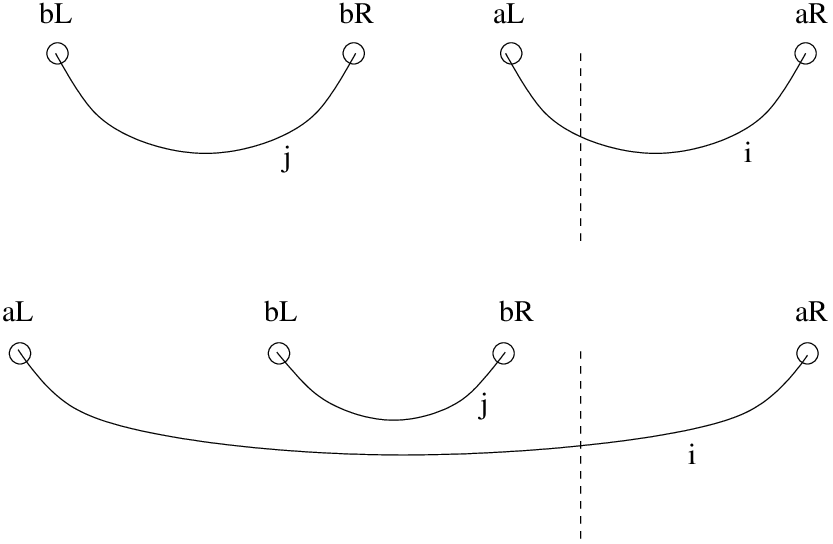}
\]
\caption{Case II of Lemma~\ref{lem:lefthandends}: Diagrams $D$ (top) and $\TLproj(D)$.}
\label{fig:lefthandendscaseII}
\end{figure}

The vertex with right-standard label $bR$ in $D$ must be to the left of
$x$, since this is so in $\TLproj(D)$ (again using Lemma~\ref{lem:recurs1}),
so, since $y$ is right-standard labelled $aL$ in $D$, the vertex
with right-standard label $bR$ in $D$ must be to the left of the
vertex with right-standard label $bL$ in $D$.
Let $i$ (respectively, $j$) be the secondary label of the arc with end-points
$aL$ and $aR$ (respectively, $bR$ and $bL$) in $D$.
By Lemma~\ref{lem:incomparable}, $i$ and $j$ are not descendants of each
other in $\tau_R(D)$. But this contradicts the fact that the
right-corresponding arcs (also labelled $i$ and $j$) in $\TLproj(D)$ are
nested. See Figure~\ref{fig:lefthandendscaseII}. The dashed vertical
line is drawn between vertices $x$ and $y$ (so the right-agreeing part of
$D$ and $\TLproj(D)$ is to the right of this line).

Case (III): Suppose that $y$ is the right-hand end of an arc in $D$ but the
left-hand end of an arc in $\TLproj(D)$. 
Let $aL$ and $aR$ be the right-standard labels of the arc in $D$ incident
with $y$ in $\TLproj(D)$ and let $bL$ and $bR$ be the right-standard labels of
the arc incident with $y$ in $D$. Note that $aR$ must be $x$ or to its
right in $\TLproj(D)$, so the vertex with right-standard label $aR$ in
$\TLproj(D)$ must also be $x$ or to the right of $x$, since $D$ and
$\TLproj(D)$ right-agree up to vertex $x$, using Lemma~\ref{lem:recurs1}. 

Since the vertex with right-standard label $bR$ in $D$ is to the left of
$x$, its right-correspondent in $\TLproj(D)$ must also be to the left of $x$
(using Lemma~\ref{lem:recurs1}) and thus the whole of the arc with right-standard labels $bL$ and $bR$ must be to the left of the vertex with right-standard
label $aL$ in $\TLproj(D)$.
Let $i$ (respectively, $j$) be the secondary label of the arc with end-points
$aL$ and $aR$ (respectively, $bR$ and $bL$) in $D$. These are the left-standard
labels of the right-corresponding arcs in $\TLproj(D)$, and by the above neither
is a descendant of the other in $\tau_R(D)$ by the definition of $\TLproj(D)$.

We claim that, using the definition of $\tau_R(D)$, $j$ is a descendant of
$i$ in $\tau_R(D)$: a contradiction. Since $\TLproj(D)$ and $D$ right-agree up to
vertex $x$, the diagram for $D$ can be drawn with no crossings to the right
of a vertical line $V$ drawn between vertices $x$ and $y$.

Let $a_1R,a_2R,\ldots ,a_kR$ be the right-hand end-points of the
arcs of $D$ with right-hand end-point at $x$ or to its right and
left-hand end-point to the left of $x$, with $a_1<a_2<\cdots <a_k$.
Note that the left-hand end-point of each of these arcs is to the left
of the vertex with right-standard label $bR$ in $D$.
It follows that arc $j$ is in the subdiagram $D^{a_1}$.

\begin{figure}
\[
\includegraphics[width=3in]{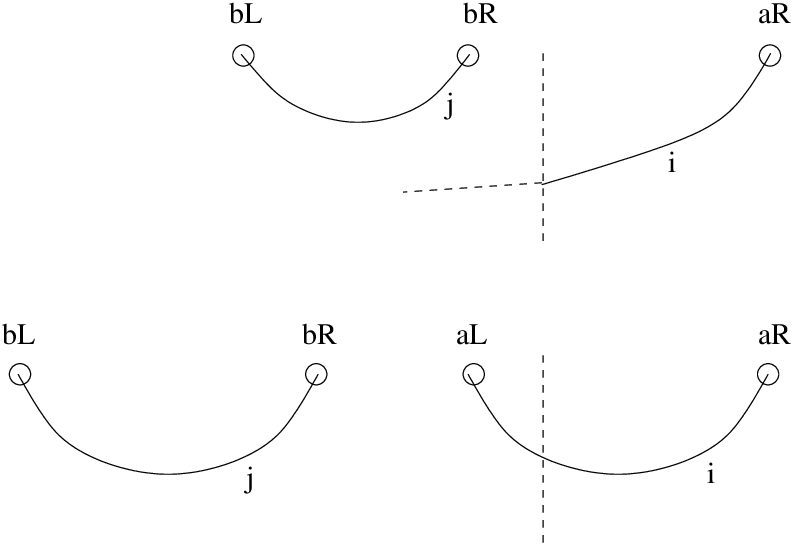}
\]
\caption{Case III of Lemma~\ref{lem:lefthandends}: Diagrams $D$ (top) and $\TLproj(D)$.}
\label{fig:lefthandendscaseIII}
\end{figure}

It follows similarly that arc $j$ is in the subdiagram $D^{a_2}$, and
by continuing to argue in this way we eventually obtain that arc $j$ is
in the subdiagram $D^i$ and thus is a descendant of $i$ as required.
See Figure~\ref{fig:lefthandendscaseIII}.

We have thus ruled out all other possible configurations and can conclude
that the lemma holds.
\Qed
%}}}

%{{{ lemma not in a valley
%\begin{lemnum}
\p \textbf{Lemma.} \label{lem:novalley}
\emph{Let $D\in J_n$, and suppose that $D$ and $\TLproj(D)$ maximally right-agree up
to a vertex $x$. Let $y$ be the vertex immediately to the left of $x$
and let $z$ be the vertex immediately to the left of $y$. Then
$z$ is the left hand end-point of an arc in $\TLproj(D)$.}
%\end{lemnum}

\emph{Proof:}
For a contradiction, we suppose that the vertex $z$ in $\TLproj(D)$ is
the right-hand end of an arc. Let its right-standard label be $cR$.
The vertex which is right-standard labelled $cR$ in $D$ must occur
to the left of vertex $x$ in $D$, as it does in $\TLproj(D)$. 
Let $dL$ be the right-standard label of the vertex $y$ in $D$
(note that by Lemma~\ref{lem:lefthandends} this vertex must be the
left-hand end-point of an arc in $D$).

We see that:

(*) The vertex with right-standard label $cR$ in $D$ occurs to the left
of the vertex with right-standard label $dL$.

The right-hand end of the arc incident with $dL$ must be either $x$ or
to the right of $x$ in $D$. This vertex is labelled $dR$.
Then the vertex with right-standard label $dL$ must be to the left
of $x$ in $\TLproj(D)$ and the vertex with right-standard label $dR$ must
be $x$ or to the right of $x$ in $\TLproj(D)$ (as both of these hold for
the right-corresponding vertices in $D$).

Since $\TLproj(D)$ has no crossings of arcs, the vertex with right-standard
label $dR$ in $\TLproj(D)$ must be to the right of the vertex right-standard
labelled $bR$ in $\TLproj(D)$ and the vertex with right-standard label
$dL$ in $\TLproj(D)$ must be to the left of the vertex with right-standard
label $cL$ in $\TLproj(D)$. Thus

(**) In $\TLproj(D)$, the arc with end-points $dL$ and $dR$ contains the
arc with end-points $cL$ and $cR$.

Let $i$ (respectively, $j$) be the secondary label of the arc with end-points
$bL$ and $bR$ (respectively, $cL$ and $bR$) in $D$. Let $k$ be the
secondary label of the arc with end-points $dL$ and $dR$ in $D$; this
coincides with the left-standard label of the arc in $\TLproj(D)$ with
these end-points, by the definition of secondary label. Then
by (**) above and the definition of $\TLproj(D)$, $j$ is a descendant
of $k$ in $\tau_R(D)$. But by (*) and Lemma~\ref{lem:incomparable},
$j$ is not a descendant of $k$ in $\tau_R(D)$, a contradiction.

It follows that vertex $z$ must be the left hand end-point of an arc in
$\TLproj(D)$ as required.
\Qed
%}}} 

%{{{ lemma

\section{Main Result}
\label{s:mainresult}
\p
In this section we will define the map $\Phi:J_n\rightarrow \Y_n$ and show that
it is an inverse to the tiling map $\Psi$, thus proving our main result,
that there is a bijection between overhang paths and Brauer diagrams.
We first need a key lemma:

%\begin{lemnum}
\p \textbf{Lemma.} \label{lem:recurs2}
\emph{Let $D \in J_n$, and suppose that $D$ and $\TLproj(D)$ maximally right-agree up
to vertex $x$. Let $r$ be as in Lemma~\ref{lem:recurs1}. Then: \\
(a) The right hand ends of the arcs in $D$ and $\TLproj(D)$ with
secondary label $r+1$ lie in the right-agreeing right-hand-end of the
two diagrams.
(b) The left-hand end of the arc with secondary label $r+1$ in $D$
is \emph{further} from the right-hand end of $D$ than the left-hand
end of the arc with the same left-standard (i.e.\ secondary) label in
$\TLproj(D)$.}
%\end{lemnum}

{\em Proof:}
The arc whose left-hand end-point is immediately to the left of vertex
$x$ in $\TLproj(D)$ (see Lemma~\ref{lem:lefthandends}) has left-standard
(i.e.\ secondary)
label $r+1$ by definition of left-standard labelling. Therefore its right-hand
end-point is in the right-agreeing part of $D$ and $\TLproj(D)$. The arc in
$\TLproj(D)$ with secondary label $r+1$ is the right-corresponding arc in
$D$. Since $D$ and $\TLproj(D)$ right-agree up to $x$, its right-hand end
must also lie in the right-agreeing part, and (a) is shown.
Since $D$ and $\TLproj(D)$ maximally agree up to vertex $x$, the left-hand
end of the arc with secondary label $r+1$ in $D$ cannot be the
vertex immediately to the left of vertex $x$, but it cannot be in the
right-agreeing part of $D$ and $\TLproj(D)$ since the left-hand end of
the right-corresponding arc in $\TLproj(D)$ does not lie in this right-agreeing
part. The result follows.
\Qed

\p
Let $X_D$ denote the number of steps to the right that the left
hand end of arc with secondary label $r+1$ in $D$ in Lemma~\ref{lem:recurs2}
above would have to be moved in order to right-agree with that arc in
$\TLproj(D)$. Define $\delta D$ to be the diagram resulting from moving the
left-hand end of the arc with secondary label $r+1$ in this way.

\p
For example, if $D$ is the diagram in Figure~\ref{fig:phi'tl}(a), then
$r=0$ and $D$ and $\TLproj(D)$ maximally right agree up to the third
vertex from the right. Since moving the left hand end point of the
arc with secondary label $1$ in $D$ five steps to the right would make it
right agree with the arc with the same label in $\TLproj(D)$, we have $X_D=5$.

\p
This means that $\delta D$ and $\TLproj(D)$ exhibit greater right-agreement,
that is, right-agreement up to a vertex to the left of the right-agreeing
part of $D$ and $\TLproj(D)$. (Note that in the example, $\delta D$ and
$\TLproj(D)$ maximally right-agree on the nine rightmost vertices).
We define $D(r) = D$ and $X_{r}=X_{D}$. Next we define
$D(r')=\delta D$ and $X_{r'}=X_{\delta_D}$
(where $r'$ is as in Lemma~\ref{lem:recurs1} for
$\delta D$), and so on, iterating the procedure.
We thus obtain a sequence
\[
D \stackrel{X_{D(r)} = X_D}{\rightarrow} 
D(r) \rightarrow D(r') \rightarrow ...
\rightarrow \TLproj(D).
\]
For any $j$ not appearing in this sequence (i.e.\ no adjustment is
required to bring the arcs with secondary label $j+1$ into right-agreement),
we define $D(j)$ to be $D(s)$, where $s$ is minimal with the property
that the arcs with secondary label $j+1$ agree in $D(s)$ and $\TLproj(D)$.
In such cases, we set $X_j=0$.

\p
Each single step counted by $X_D=X_r$ can be implemented on $D$ by an
adjacent pair permutation of vertices. Thus by extending the diagram
above to include a single crossing $\sigma_{j(r)}$ (say) in the appropriate
position, for each such step, we can build up the transformation 
$D \rightarrow \delta D$. We repeat this procedure for each
transformation in the above sequence.

\begin{figure}
\[
\includegraphics[width=3.3in]{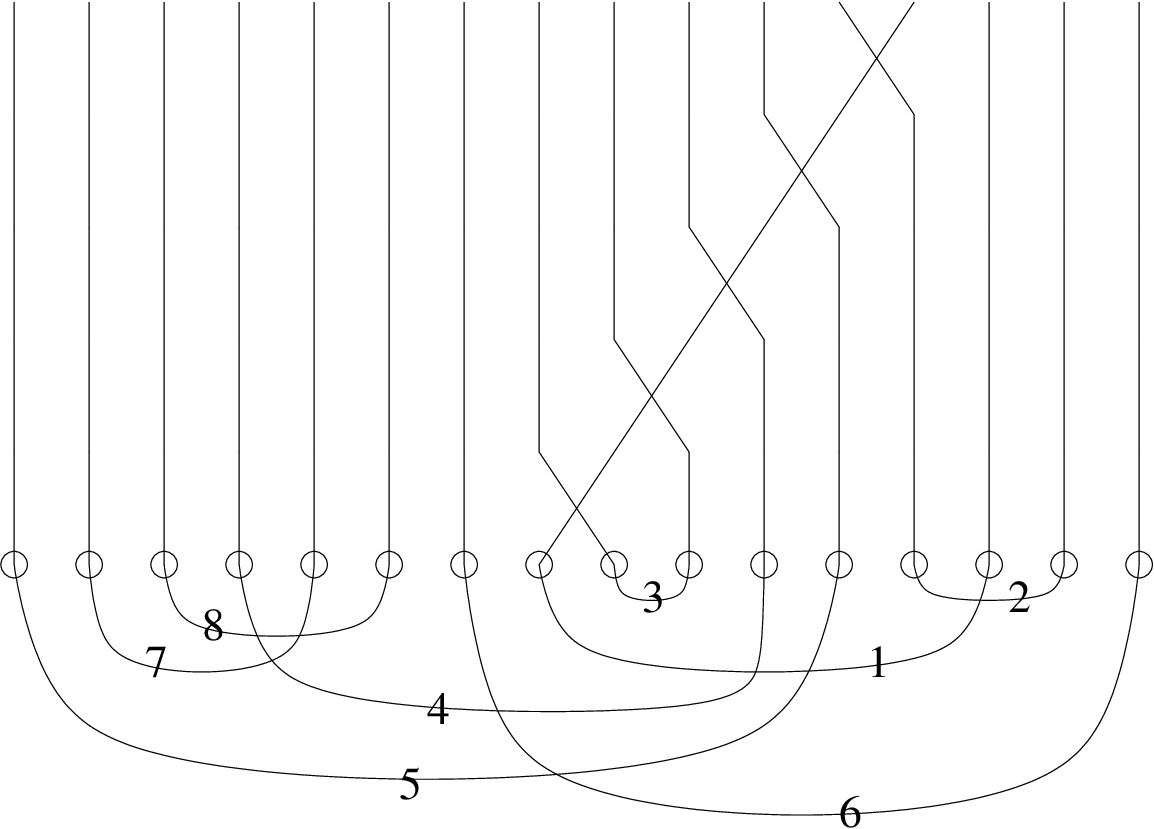}
\]
\caption{The first part of the extension for $D$ as in Figure~\ref{fig:phi'tl}(a).}
\label{fig:extension}
\end{figure}

In our example, the extension for the first transformation is shown in
Figure~\ref{fig:extension}.

\p
Of course  $\sigma_{j(i)}^2=1$,
so applying the collection of these changes in reverse order to
$\TLproj(D)$ brings us to $D$. 

\p \defcom{The inverse map, $\Phi$}
\label{def:inversemap}
Given a Brauer diagram, $D$, let $\Phi(D)$ be the diagram obtained by
starting from the \Dyck\ path $\TLinv(\TLproj(D))$ for $\TLproj(D)$
and, for each $i$, appending to the step in the Dyck path of form $(1)$
with label $i$  a left-overhanging stack of tiles of length $X_{D(i)}$.
We remark that it follows from Lemma~\ref{lem:novalley} that $\Phi(D)$
is an \overhang\ path.

\p
It follows from the above that applying $\Psi$ to $\Phi(D)$ we obtain
the original diagram $D$. (An example follows shortly).

\p
It is well known (see also Section~\ref{s:easymap})
that the cardinality of $J_n$ is equal to $(2n-1)!!$
so it follows from Lemma~\ref{lem:grimecardinality} that $|J_n|=|\Y_n|$.
We have therefore shown that:

\begin{theo*}\label{th:main}
The maps $\Phi$ and $\Psi$ are inverse bijections between the set $J_n$ of
Brauer diagrams of degree $n$ and the set $\Y_n$ of \overhang\ paths of degree $n$.
\end{theo*}

%}}}
%{{{ example improving partial agreement

\section{Example of a $\delta$-sequence.} 
\label{s:example}
\p
We now give an example demonstrating the main theorem.

\begin{figure}
\[
\includegraphics[width=4in]{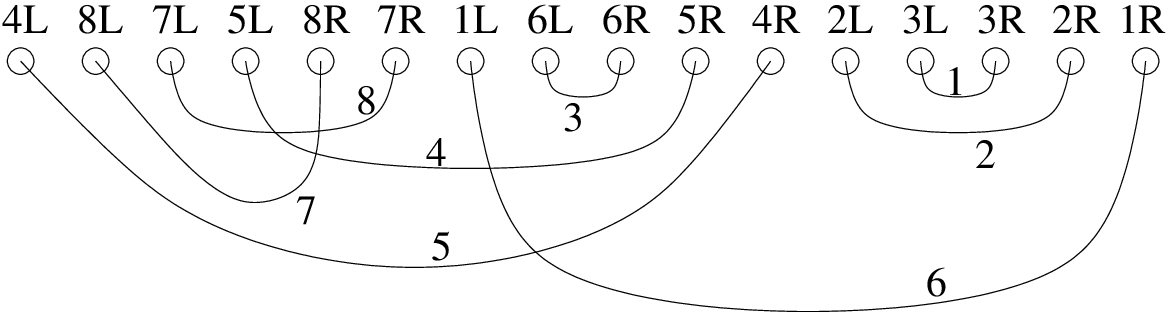}
\]
\caption{The diagram $D(3)=\delta D$ for $D$ as in Figure~\ref{fig:phi'tl}(a).}
\label{fig:D3}
\end{figure}

\p
Let $D$ be the Brauer diagram in Figure~\ref{fig:phi'tl}(a). We have seen
already that $\TLproj(D)$ is the TL diagram in Figure~\ref{fig:phi'tl}(b),
and that $D$ and $\TLproj(D)$ maximally right agree up to the
third vertex from the right, including $r=0$ arcs in their entirety.
Thus we write $D(0)=D$. We have seen that $X_D=X_0=5$, so we first move
the arc secondary labelled $1$ five steps to the right to obtain diagram
$\delta D$; see Figure~\ref{fig:D3}.

\begin{figure}
\[
\includegraphics[width=4in]{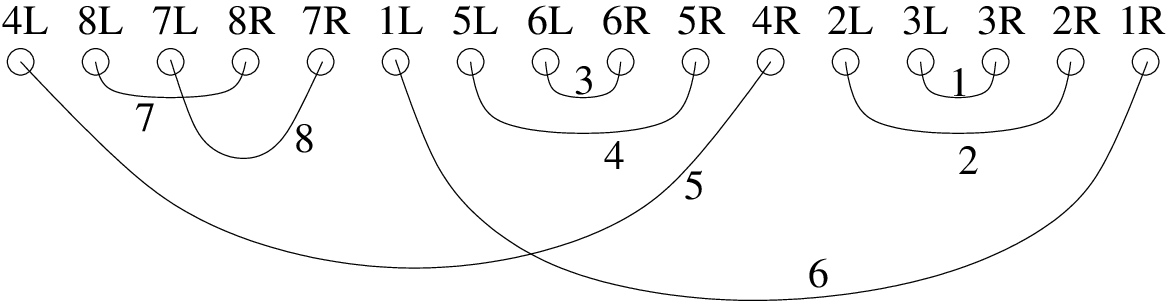}
\]
\caption{The diagram $\delta D(3)=D(4)$}
\label{fig:D4}
\end{figure}

\p
We observe that $\delta D$ and $\TLproj(D)$ maximally agree up to the
ninth vertex from the right, including the arcs secondary-labelled $1$.
In fact the arcs secondary-labelled $2$ and $3$ also right-agree, so
$r'=3$ and we write $\delta D=D(3)$ (thus $D(1)=D(2)=\delta D$).
By moving the arc secondary-labelled $r'+1=4$ three steps to the right we can
make it right-agree with the arc with the same secondary label in
$\TLproj(D)$, giving the diagram $\delta D(3)$ shown in Figure~\ref{fig:D4}.
Thus $X_3=3$.

\begin{figure}
\[
\includegraphics[width=4in]{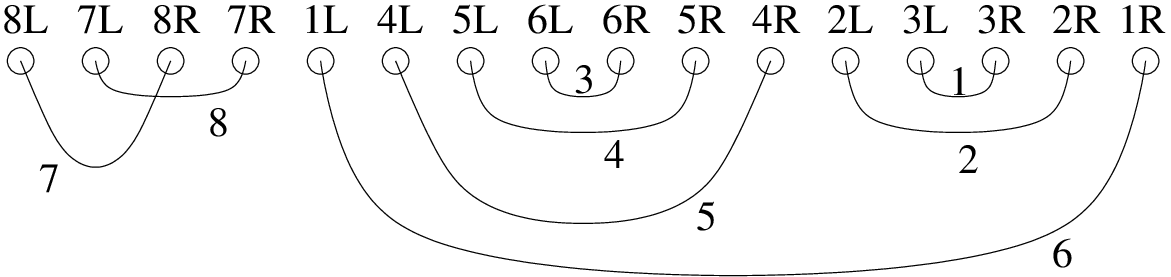}
\]
\caption{The diagram $\delta D(4)=D(6)$}
\label{fig:D6}
\end{figure}

\p
We see that the arcs secondary-labelled $1,2,3,4$ lie in the right-agreeing
parts of $\delta D(3)$ and $\TLproj(D)$, so $\delta D(3)=D(4)$. We
compute that $X_4=5$ and diagram $\delta D(4)$ is shown in Figure~\ref{fig:D6}.

\p
Next, the arcs secondary-labelled up to $6$ lie in the right-agreeing
parts of $\delta D(4)$ and $\TLproj(D)$, so $\delta D(4)=D(6)$. We
compute that $X_6=1$ and $\delta D(6)=\TLproj(D)$.

\begin{figure}
\[
\includegraphics[width=4in]{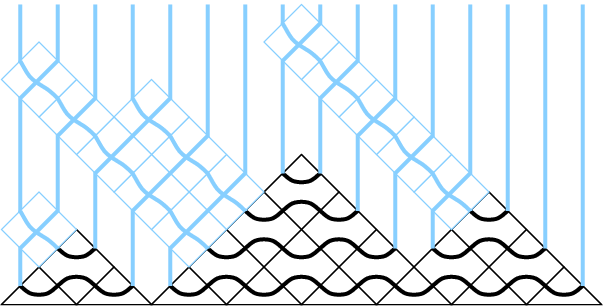}
\]
\caption{The \overhang\ path $\Phi(D)$ and its tiling.}
\label{fig:answerpath}
\end{figure}

\p
Using this data to construct $\Phi D$ we obtain the \overhang\ path in
Figure~\ref{fig:answerpath}. Note that the tiling of this path
does indeed return $D$. 

%}}}
%}}}

%{{{

\section{A simple bijection between \overhang\ paths and Brauer diagrams}
\label{s:easymap}
%{{{ 1.1

\p
In this section we give a simple bijection between \overhang\ paths and
Brauer diagrams, also given by tiling. We note, however, that it does not have
the property that it restricts to the tiling bijection for the Temperley-Lieb
case, described in Section~\ref{s:TLtile}.

Let $n\in\mathbb{N}$. Recall that $J_{n-1}^1$ denotes the set of partial
Brauer diagrams with $n-1$ pairs and one singleton. There is a bijection
$$s_{2n}:J_n \rightarrow J_{n-1}^1,$$
given by deleting the rightmost vertex of a Brauer diagram.
The inverse adds a single vertex at the right hand end and joins it
with the singleton.
There is a map from $J_{n-1}^1$ to $J_{n-1}$
obtained by deleting the singleton. There are $2n-1$ possibilities for
the singleton, giving a bijection:
$$s_-:J_{n-1}^1 \rightarrow J_{n-1} \times \ul{2n-1},$$
where $\ul{k}$ denotes the set $\{1,2,\ldots ,k\}$ for any $k$.
Thus $d_{2n}:=s_{-}\circ s_{2n}$ is a bijection from $J_n$ to
$J_{n-1}\times \ul{2n-1}$. It follows that
$$|J_n| = (2n-1) |J_{n-1}|,$$
and thus that $|J_n|=(2n-1)!!$.
%}}}
%{{{ 3
It follows from the above that
$$\kappa:=d_2\circ d_4 \circ \cdots \circ d_{2n}:J_n \rightarrow A_n,$$
where
$$A_n := \{ (x_1,x_2,...,x_n) \in \Z^n \; | \; 1\leq x_i < 2i \}$$
is a bijection.
%}}}
%{{{ Dyck

\p
In this section only, we shall regard a \emph{Dyck path}
as a walk on $\Z\times \Z$ from $(0,0)$ using
steps from $\{ (1,0),(0,1) \}$, such that the walk never drops below the
line parallel to the vector $(1,1)$
(equivalently, if the {\em height} of a point $(x,y)$ is defined to be $y-x$,
negative heights are not allowed).
It is clear that such a path can be
transformed into a Dyck path as defined in Section~\ref{s:defnsetc}
by rotating it through $45$
degrees clockwise about the origin and stretching it by a factor of
$\sqrt{2}$. We consider such paths whose end-point is $(n,n)$.

%}}}
%{{{ overhang

\p
Similarly, in this section only, we shall regard an \emph{\overhang\ path}
as a generalisation of such a walk in which steps of
the form $(-1,0)$ are also allowed, but the walk also never drops below
(i.e.\ to the left of) the line defined by the $(-1,1)$ vector (and
the path may not visit the same vertex twice).
Such a path is characterised by the sequence of $x$-coordinates of its
$(0,1)$-steps. 
The first entry in this sequence is necessarily $0$, the second
lies in $\{ -1,0,1\}$, the third lies in $\{-2, -1,0,1,2 \}$, and so on. 
(In the Dyck path case the negative positions do not occur.)
Let $O_n$ denote the set of such paths ending at $(n,n)$.

%}}}
%{{{ bij

\p
It is clear that there is a bijection from $A_n$ to $O_n$ taking
an element
$(x_1,x_2,\ldots ,x_n)$ of $A_n$ to the \overhang\ path with sequence
of $x$-coordinates of its $(0,1)$-steps given by $x_i-i$, $i=1,2,\ldots ,n$.

\p
We have thus constructed a bijection
\[
J_n \rightarrow A_n \rightarrow O_n
\]

%}}}
%{{{ inv bij
One way to construct the inverse of the above bijection is to start with an
element of $O_n$ and to regard this as a partial tiling of the plane with 
$1 \times 1$  tiles.  That is, one fills the interval between a given
\overhang\ path and the lowest path with tiles. 
One also tiles the interval between the $(1,1)$ line
and the lowest path with half-tiles in the obvious way. 
One then decorates all the square tiles with crossed
lines from edge to opposite edge; and the triangular tiles each
with a single line from short edge to short edge. 
This decoration gives the corresponding element of $J_n$. 

%}}}

\section{The Temperley-Lieb/Dyck path paradigm} \label{s:TLparadigm}

%{{{ local macros
\newcommand{\gen}{g}  %% braid gen
\newcommand{\DD}{{\mathsf D}} 		%%% ariki-koike-levy algebra
\def\Cy{{C}}				%%% cyclic group
\newcommand{\yy}{y}			%%% $q$-symmetriser
\def\expchoo(#1,#2){{\tiny \left( \!\!\!\!\! \begin{array}{c} #1 \\ #2 \end{array} \!\!\!\!\! \right) }}
\def\fracalt#1#2{\left( \! #2 \! \right)^{-1} #1} 	%%%denominator as an inverse
\def\Specht#1{{\cal S}^{#1}}
\def\eps{\epsilon}			%%% g_1 eigenvalues
\def\Bred{B^{\mbox{{\tiny red}}}}	%%% reduced word basis
%\newcommand{\QED}{$\Box$}
%}}}

%{{{ Representations \label{irrep0}
\label{s:orthogonalform}
\p
In this section we explain how Dyck paths and related walks arise in the
representation theory of the Hecke algebra, via representations
arising from outer product representations of the symmetric group.

%{{{ preamble
\p
Fix $q\in \mathbb{C}$ and $n\in\mathbb{N}$.
Let $H_n = H_n(q)$ denote the usual Hecke
algebra of degree $n$ over $\mathbb{C}$ (we work over $\mathbb{C}$ for
simplicity). Thus $H_n$ is the $\mathbb{C}$-algebra with generators
$g_1,g_2,\ldots g_{n-1}$ subject to relations $g_ig_j=g_jg_i$ if
$|i-j|>1$, $g_ig_jg_i=g_jg_ig_j$ if $|i-j|=1$, and $(g_i+q^2)(g_i-1)=0$.

%}}}
%{{{ de:multipartn
\p
\defn
For $d\in \mathbb{N}$, we denote by $\Gamma^d_n$ the set of all $d$-tuples
$\lambda=(\lambda^1 ,\lambda^2 ,...,\lambda^d )$ of Young diagrams,
with $|\lambda|= \sum_{i=1}^d |\lambda^i |=n$.  
%

%}}}
%{{{ We now review 

Fix $\lambda\in \Gamma^d_n$. Then a {\em tableau} of shape $\lambda$ is
any arrangement of the `symbols' $1,2,..,n$ in the $n$ boxes of $\lambda$.
Such a tableau is said to be \emph{standard} if each component tableau
is standard. We denote the set of all standard tableaux of shape
$\lambda$ by $T^{\lambda}$.

\p
We number the rows of $\lambda$ by placing the whole of the component diagram
$\lambda^{i+1}$ under $\lambda^i$ for all $i$, and numbering the
rows from top to bottom.
We then define a total order $<$ on standard tableaux of shape $\lambda$
by setting $T<U$ if the highest number which appears in different rows of
$T$ and $U$ is in an earlier row in $U$.

%}}}
%{{{ Let $T^{\lambda}_p$ be
\p
Let $T$ be a tableau. For $i\in \{1,2,\ldots ,n-1\}$, let
$\sigma_i=(i\ i+1)\in \Sigma_n$, the symmetric group of degree $n$.
We define $\sigma_i (T)$ to be the tableau obtained by
interchanging $i$ and $i+1$. In this way we get an action of $\Sigma_n$
on the set of all tableaux of shape $\lambda$, but we note that
this action does not necessarily take a standard tableau to a
standard tableau.

%}}}
%{{{ prop p1q 

\p
Let $x=(x_1,x_2,\ldots ,x_d)\in \mathbb{R}^d$. For $i,j\in\{1,2,\ldots ,n-1\}$
and $T\in \Gamma^d_n$ let $h^x_{ij}=h^x_{ij}(T)$
denote the {\em generalised hook length} between the symbols $i$ and $j$
in $T$. Thus $h^x_{ij}$ is given by:
\[
h^x_{ij} = h^0_{ij} +x_{\# i} -x_{\# j},
\]
where $h^0_{ij}$ is the usual hook length obtained by superimposing the 
component tableaux of $T$ containing $i$ and $j$,
and $\# i$ is the number of the component containing $i$ in $T$.
See~\cite{MartinWoodcockLevy2000}
(note that there is a typographical error in this paper at the relevant point)
and also~\cite[p.244]{Martin91}. Geometrically, one may think of putting all
the individual Young diagrams $\lambda^i$ in the same plane, each with its top
left box in position $(0,x_i)$.

For an integer $m$, we will write, as usual,
$$[m]=\frac{q^m-q^{-m}}{q-q^{-1}}.$$

%\begin{prnum}
\p \textbf{Proposition}
\label{prop:outerproduct} \label{pa:10.6}
 {\rm \cite{MartinWoodcockLevy2000}}.
\emph{
Let $\lambda\in \Gamma^d_n$ and assume that $[h_{i,i+1}^x(T)]\not=0$ for all
$i\in \{1,2,\ldots ,n-1\}$ and all $T\in T^{\lambda}$.
Then the set $T^{\lambda}$ is a basis for a left $H_n$--module $R^{\lambda}$.
For $i\in\{1,2,\ldots ,n-1\}$ and $T\in T^{\lambda}$, the action is as follows: \\
(a) If $i,i+1$ lie in the same row of $T$ then $\gen_i T = T$. \\
(b) If $i,i+1$ lie in the same column of $T$ then $\gen_i T = -q^2T$. \\
(c) If neither (a) nor (b) hold, then $\sigma_i(T)$ is also standard.
Let $h=h_{i,i+1}^x$. Then the action is given by
$$\gen_i \mat{c} T^{\lambda}_p \\ \sigma_i(T^{\lambda}_p ) \tam
   = 
  \left( \mat{cc} 1 & 0 \\ 0 & 1 \tam - \frac{q}{[h]}    
     \mat{cc} {[h+1]} & {[h-1]} \\
              {[h+1]} & {[h-1]} \tam  \right)
		\mat{c} T^{\lambda}_p \\ \sigma_i(T^{\lambda}_p) \tam,
$$
provided $T<\sigma_i(T)$.}

%\end{prnum}
%
(Young's orthogonal form (see e.g.~\cite[IV.6]{boerner70})
involves an action via symmetric matrices related to those above via
conjugation).

%}}}
%{{{ new

\p \label{p:regardaswalks}
We now restrict attention to the situation in which $\lambda$ has exactly
two components, each consisting of exactly one row.
We can represent $T\in T^{\lambda}$ by an $n$-tuple $(a_1,a_2,\ldots ,a_n)$
with entries in $\{1,2\}$, defined by the condition that
$i\in \lambda^{a_i}$ for all $i\in\{1,2,\ldots ,n\}$. Such a tuple can
be regarded as a walk of length $n$ in $\mathbb{Z}^2$ starting at the origin.
The $i$th step of the walk consists of adding the vector $(1,1)$ if
$a_i=1$ or adding the vector $(1,-1)$ if $a_i=2$.

\begin{figure}
\[
\includegraphics[width=5in]{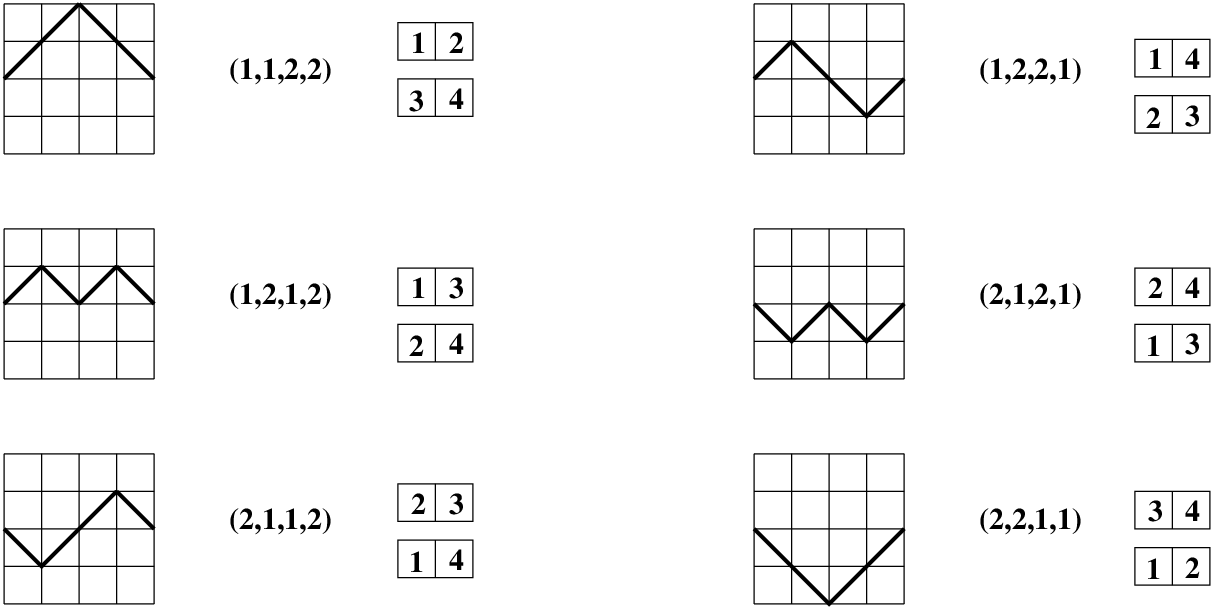}
\]
\caption{Standard tableaux of shape $((2),(2))$ and the corresponding walks.}
\label{fig:examplewalks}
\end{figure}

\p
For example, if $n=4$ and each component of $\lambda$ is a row of length
$2$, the elements of $T^{\lambda}$ and the corresponding tuples and
walks are as shown in Figure~\ref{fig:examplewalks}.

\p
We note that in the walk realisation of a standard tableau $T$,
$\sigma_i$ swaps a pair of
steps $(1,2)$ with the pair $(2,1)$, i.e.\ a local maximum is
swapped with a local minimum or vice versa.
Thus, in order for there to be mixing between two basis elements as in
Proposition~\ref{prop:outerproduct}(c), the corresponding walks
must agree in all but their $i$th and $(i+1)$st steps,
and in each diagram separately the second coordinate (or \emph{height})
after $i-1$ steps and after $i+1$ steps must coincide.

\begin{figure}
\[
\includegraphics[width=3in]{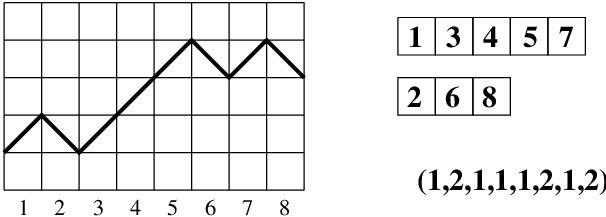}
\]
\caption{Example of a walk and corresponding standard tableau.}
\label{fig:longexamplewalk}
\end{figure}

\p
In fact, in this case the height after $i-1$ steps coincides with the
usual hook length $h^0_{i,i+1}$: $i+1$ appears in the first component of $T$
(in the $k$th box, say) and $i$ in the second component of $T$
(in the $l$th box, say) and the height of the walk after the $(i-1)$st
step is $(k-1)-(l-1)=k-l$, i.e. the hook length.
For an example, see
Figure~\ref{fig:longexamplewalk}. Here it can be seen that the height
of the walk after the $5$th step is equal to the hook length $h^0_{6,7}=3$.
In general, we have that if
$T<\sigma_i(T)$ then $h^x_{i,i+1}$ is equal to the sum of $x_1-x_2$ and the
height of the walk after $i-1$ steps.

\p
If $h^x_{i,i+1}=1$, it follows from the description of the action
in case (c) that the elements are not actually mixed.
It follows that, if we choose $x$ so that $x_1-x_2=1$, there is an action
of $H_n$ on the set of (standard tableaux corresponding to) walks which
do not go below the horizontal axis given by the formulas in
Proposition~\ref{prop:outerproduct}. In fact, in this case, the action
cannot be extended to the whole of $T^{\lambda}$ since the action is
not defined for hook length zero.

\p
Similarly, if we set $x_1-x_2=2$, only the walk $(2,2,1,1)$ is
decoupled from the rest. In other words, changing the value of $x$ allows us to
define a module for $H_n$ with basis elements corresponding to walks
which do not go below a certain "exclusion" line.

\p \label{p:Zgraph}
Let $\Gamma_{\Z}$ be the graph with vertices $\mathbb{Z}$ and edges
joining integers with difference $1$. Then 
the walks we have been considering can be regarded as walks on
$\Gamma_\mathbb{Z}$ by projecting onto the second coordinate.
Thus, in summary, we have extracted an $H_n$-module with a basis of walks
on $\Gamma_{\Z}$ which only visit vertices on a certain
subgraph, from the formal closure of a Zariski-open set of modules (that is,
actions depending on a parameter) whose bases consist of walks on a larger
subgraph. The decoupling of the subgraph, in this sense, is determined
by the structure of the graph.

\p \label{p:walkbasis}
The case $x_1-x_2=1$ is special in that the decoupled module is
irreducible for generic values of $q$. It is an analogue of setting
the usual $\rho$-shifted position of the boundary
of the dominant region in the Weyl group construction in Lie theory.
The generic irreducible module is denoted by
$\Delta_n^{TL}(\lambda^1,\lambda^2)$.
The most interesting step, however, is the next one. We now fix $x_1-x_2=1$,
and also specialise $q$ to be an $l$th root of unity, so that $[l]=0$.
In this situation, there is a further decoupling: we obtain a module
whose basis corresponds to walks whose height is bounded above by
$l-1$. In other words, we now only include walks that lie between two
`walls': the lines given by setting the second coordinate to $0$ and $l-1$.
It can be shown that this module is simple in this specialisation.
Such simple modules are otherwise very hard to extract, but here their
combinatorics is manifested relatively simply.

\p
It is this feature that we aim, eventually, to duplicate for the Brauer
algebra. Although we reiterate that in this article we have not addressed the
representation theory of the Brauer algebra --- as a first step, we
have considered a Brauer analogue of the underlying combinatorial
correspondence.

%}}}
%}}}

%%%%%%%%%%%%%%%%%% NEW BITS BELOW HERE %%%%%%%%%%%%%%%

\newcommand{\E}{{\mathbb E}}
\newcommand{\ZN}{\Z^\N}
\newcommand{\ZNp}{\Z^\N_{\rho}}
\newcommand{\ZNd}{\Z^\N_{\rho - \delta}}
\newcommand{\RN}{\R^\N}
\newcommand{\Zgraph}[1]{\Gamma_{\Z^{#1}}}
\newcommand{\Aregular}{$A$-regular}
\newcommand{\QQ}{\delta}
\newcommand{\ww}{\omega}
\newcommand{\YG}{{\mathcal Y}}
\newcommand{\ee}{e}
\newcommand{\footnot}[1]{}

%{{{ chat

\p
The evidence for an analogy with the Brauer algebra is captured by the
following summary of Cox-De~Visscher-Martin's 
geometrical machinery for the representation theory of
the Brauer algebra \cite{cdm09,cdm10,martin09}.%

\p
Let $B_n(\QQ)$ denote the Brauer algebra of rank $n$ with parameter
$\delta\in\R$.
For this case we may generalise our graph
$\Gamma_{\Z}$ (as in \ref{p:Zgraph})
and its role in representation theory as follows.
Let $\E^n = \R^n$ be Euclidean $n$-space, and $\Z^n$ the integral
lattice. Let $\E^{n} \hookrightarrow \E^{n+1}$ be the natural
inclusion and $\E^f$ denote the inverse limit, with basis $\{ e_1,
e_2,... \}$. 
%
%}}}
%{{{ standards
%
%{{{ walk
%
%Let $\E^n = \R^n$ 
%be Euclidean $n$-space and $\Z^n$ the integral lattice in
%$\E^n$. Let $\E^n \hookrightarrow \E^{n+1}$ be the natural inclusion,
%and let 
%$\E^f$ denote the inverse limit of these inclusions,
%with basis $\{ e_1, e_2, ... \}$ (and $\Z^f$ similarly).
%
%\noindent
Thus $\Z^f \subset \ZN$.
%
%\mdef We d
Define $\Zgraph{n}$ 
(for $n \in \N$, or $n=\N$)
to be the graph with vertex set $\Z^n \cup (\Z+1/2)^n$ and 
an edge $(x,x')$ if $x-x' = \pm e_i$.

%Note that  $\Zgraph{n}$ embeds naturally in $\R^n$, where the edges
%may be embedded as straight lines of length 1.

%}}}
%{{{ reflect

%\mdef
%Let $(ij)$ denote the reflection in the hyperplane which takes 
%$e_i \leftrightarrow e_j$ and fixes all other basis elements
%(and $(ij)_-$ takes $e_i \leftrightarrow -e_j$ similarly).
%Extend this action to an action on $\R^\N$ in the obvious way.

%\mdef For given $n$, we d
\p
Define $(ij)_{\pm}:\R^n \rightarrow \R^n$ by 
\[
(ij)_{\pm} (x_1, x_2, ... , x_i ,...,x_j ,...) 
           =  (x_1, x_2, ... , \pm x_j ,...,\pm x_i ,...)
\] 
%and define $(ij)=(ij)_+$.
Write $H_{(ij)_\pm}$ for the reflection hyperplane in $\R^n$ associated
to $(ij)_\pm$;
$\H_A = \{ H_{(ij)_+} \}_{ij}$; 
and $\H = \{ H_{(ij)_{\pm}} \}_{ij}$. 
The open (codimension 0) components  
of $\RN\setminus \H_A$  are  {\em chambers}.

\p
By definition, a {\em regular part} of $\Zgraph{n}$ is
the full subgraph on vertices lying in a fixed 
chamber.
A vertex in $\Zgraph{n}$ is {\em \Aregular} if it lies on no hyperplane
of form $H_{(ij)_+}$. 
A walk on  $\Zgraph{n}$ is {\em \Aregular} if it visits 
only \Aregular\ vertices. 

%}}}

\p
The TL module bases we have been reviewing correspond 
(after some tweaking) to the case $n=1$.
We now consider  $n=\N$. 

%{{{ N

%\mdef We d
\p
Define  % $\ww$, $\rho$ in $\R^\N$ by 
$
-2\ww=(1,1,1,...) 
$,
$
-\rho = (0,1,2,...)
$
and, for any $\delta\in\R$,
$
\rho_{\QQ} = \QQ \ww + \rho
$  
in $\R^\N$.
%\mdef 
\label{de:intseq}
%We shall call a sequence $v$ `integral' if either $v$ or $v+\ww$
%lies in $\Z^\N$.
%
\footnot{Asides:
%{{{ index

\mdef
We call $v \in \RN$ {\em rational} if 
\[
f_v(x) \; := \; \sum_i v_i x^i
\]
is the Maclaurin series of a rational function. For example
any $v \in \E^f$ is rational, and 
$\ww$ and 
$\rho$
%\[ 2\ww \; := \; -(1,1,1,...) \]
%and \[ \rho_0 \; := \; -(0,1,2,...) \]
are rational.

%}}}
%{{{ dom chamber

\mdef
The subspace of $\RN$ fixed by any 
reflection $(ij)$ partitions $\RN$ into three
parts --- two open parts and the subspace itself.
The collection of all such subspaces partitions $\RN$ into open parts
of codimension zero
(called chambers) and various facets. 
For example the set of strictly descending sequences form a chamber.
We shall call it the `dominant' chamber.
The facets bounding this chamber are various intersections of
hyperplanes of form $(i,i+1)$.

%}}}
}
%
%}}}
%{{{ rho delta A-reg
%
%\mdef
Note that $\rho_{\QQ}$ is \Aregular\ for any $\QQ$. 
We call the chamber containing $\rho$ 
%(and hence, by continuity, all the $\rho_{\QQ}$s)  
the {\em dominant chamber}.
We call the 
full subgraph of $\Zgraph{\N}$ with vertices in the dominant chamber 
the {\em dominant part} of $\Zgraph{N}$,
and denote it $(\Zgraph{\N})^{\rho}$.

%\mdef
%We shall regard the natural injection $\Lambda \hookrightarrow \Z^f$
%as an inclusion.

\p
For $\delta \in \R$ define
$
\ee_{\delta} : \Z^f \hookrightarrow \RN
$
by 
$\lambda \mapsto \lambda + \rho_{\delta}$; 
%
%Note that the images  $e_{\QQ}(\Lambda)$ 
%all lie in the dominant chamber. Note also that no two of these 
%images intersect (as $\delta$ varies); and that the image is a set of
%integral sequences 
%(in the sense of (\ref{de:intseq}))
%for every integer $\delta$. 
%
%\mdef For each $\delta \in \Z$ define the set of  $\delta$-finite weights 
%Define
$
X_{\delta} 
= \ee_\QQ (\Z^f )   %%\{ -\rho_\QQ +x \; | \; x \in  \Z^{f} \} .
$ 
%
%}}}
%{{{ connected graph etc
%
%\mdef
%We say that two vertices on a graph are connected if there is a finite
%walk (a finite sequence of edges) between them. 
%A graph is connected if every pair of vertices is connected.
%Thus for example $\YG$ is connected,
%$\Zgraph{\N}$ is not connected.
%
%\mdef For $\delta\in\Z$ we 
and
define $(\Zgraph{N})^{\rho}_{\rho_\QQ}$ to be the connected component of $(\Zgraph{N})^{\rho}$ containing
$\rho_\QQ$. 

%Note that $\ZNd{\QQ} $ intersects  $\ZNd{\QQ'} $ only if $\QQ=\QQ'$.

%}}}
%{{{ claims

%{\mlem{
%\begin{lemnum}
\p \textbf{Lemma.} \label{lem:YGiso1} %CLAIM. 
For $\QQ\in\Z$, $\ee_\QQ$ extends to a map
from 
the Young graph $\YG$ 
to $\Zgraph{\N}$,
%%defined by $\lambda \mapsto \rho_{\QQ} + \lambda$ 
inducing an isomorphism from 
$\YG$ to $(\Zgraph{N})^{\rho}_{\rho_\QQ}$. 
%the connected component containing 
%$\rho_{\QQ}$  of the dominant part of $\Zgraph{\N}$. 
%%%%%containing $\rho_{\QQ}$ 
%%%is isomorphic to $\YG$. 
%\end{lemnum}
%}}

%For given $\QQ$, we shall call the Young diagram $\lambda$ corresponding to any
%vertex $v = \ee_\QQ (\lambda)$ 
%of $\ZNd{\QQ}$ under this isomorphism the {\em preimage} of $v$.

%{\mth{ %CLAIM
\p \textbf{Theorem} \label{th:brauercase}
\cite{cdm} (Corollary). 
The set of \Aregular\ walks of length $n$ on $\Zgraph{\N}$ from
$\rho_\QQ$ to $\rho_\QQ + \lambda$ is a basis for the
$B_n(\QQ)$-module $\Delta_n(\lambda)$ (for any $\QQ$). 
%}}
%\proof{This follows directly from Lemma~\ref{lem:YGiso1} 
%and (\ref{de:Wbas1}).}

\p
{ Note that we are not quite in a position to consider
{\em unrestricted} walks as in \ref{p:regardaswalks} 
here,
%of length $n$ on $\Zgraph{\N}$ from
%$\rho_\QQ$ to $\rho_\QQ + \lambda$ 
%(the `moral equivalent' of a basis for a Young module)
since there are
are infinitely many.}

\p
Theorem~\ref{th:brauercase} is an analogue of the geometrical walk basis
construction of $\Delta^{TL}_n(\lambda^1,\lambda^2)$ in \ref{p:walkbasis}.
However here Cox-De~Visscher-Martin do not give an action
(the proof is abstract representation theory).
They go further and give bases for Brauer simple modules for fixed
$\QQ\in\Z$ in terms of walks further restricted by $H_{(ij)_-}$
hyperplanes (an analogue of \ref{p:walkbasis}),
but again an {\em action} does not follow automatically.
What is needed is an analogue of the geometrical (hook length) mixing
criterion. 

\p
Altogether then, Cox-De~Visscher-Martin give rather strong evidence that
there is a geometrical walk-based construction for simple modules,
but leave crucial pieces of the guiding
analogy with TL to be filled in.
One of these is a tiling map, and that is what we have provided here.

%}}}

%}}}

%}}}

%}}}

\end{document}